\documentclass[%
  a4paper,
  onecolumn,
]{mypreprint}

\usepackage[english]{babel}

\usepackage{graphicx}
\usepackage{subcaption}
\usepackage{tikz}

\usepackage{pgfplots}
  \pgfplotsset{compat = 1.13}
  \usepgfplotslibrary{external}
  \tikzset{external/system call = {%
    pdflatex \tikzexternalcheckshellescape
      -halt-on-error
      -interaction=batchmode
      -jobname "\image" "\texsource"}}
  \tikzexternaldisable

\usepackage{amsmath}
\usepackage{amssymb}
\usepackage{amsthm}

\usepackage{algorithm}
\usepackage{algpseudocode}

\newcommand{\trans}{\ensuremath{\mkern-1.5mu\mathsf{T}}}
\newcommand{\herm}{\ensuremath{\mathsf{H}}}

\DeclareMathOperator{\real}{Re}

\renewcommand{\vec}{\ensuremath{\mathrm{vec}}}
\newcommand{\grad}{\nabla}
\newcommand{\opt}{\ensuremath{\mathrm{opt}}}
\DeclareMathOperator{\conv}{conv}

\newcommand{\R}{\ensuremath{\mathbb{R}}}
\newcommand{\C}{\ensuremath{\mathbb{C}}}
\newcommand{\N}{\ensuremath{\mathbb{N}}}
\newcommand{\K}{\ensuremath{\mathcal{K}}}

\renewcommand{\i}{\ensuremath{\mathfrak{i}}}

\newcommand{\nG}{\ensuremath{n}}
\newcommand{\mG}{\ensuremath{m}}
\newcommand{\pG}{\ensuremath{p}}
\newcommand{\nK}{\ensuremath{n_{\mathrm{K}}}}

\newcommand{\nopt}{\ensuremath{N}}

\newcommand{\maxdisc}{\ensuremath{L}}

\newcommand{\AK}{\ensuremath{A_{\mathrm{K}}}}
\newcommand{\BK}{\ensuremath{B_{\mathrm{K}}}}
\newcommand{\CK}{\ensuremath{C_{\mathrm{K}}}}
\newcommand{\DK}{\ensuremath{D_{\mathrm{K}}}}
\newcommand{\EK}{\ensuremath{E_{\mathrm{K}}}}
\newcommand{\xK}{\ensuremath{x_{\mathrm{K}}}}
\newcommand{\dxK}{\ensuremath{\dot{x}_{\mathrm{K}}}}

\newcommand{\Aclosed}{\ensuremath{A_{\mathrm{c}}}}
\newcommand{\Bclosed}{\ensuremath{B_{\mathrm{c}}}}
\newcommand{\Cclosed}{\ensuremath{C_{\mathrm{c}}}}
\newcommand{\Dclosed}{\ensuremath{D_{\mathrm{c}}}}
\newcommand{\Eclosed}{\ensuremath{E_{\mathrm{c}}}}
\newcommand{\Gclosed}{\ensuremath{G_{\mathrm{c}}}}

\newcommand{\Hinf}{\ensuremath{\mathcal{H}_{\infty}}}
\newcommand{\Linf}{\ensuremath{\mathcal{L}_{\infty}}}

\newcommand{\hfgs}{\mbox{HF-GS}}
\newcommand{\rmfgs}{\mbox{RMF-GS}}
\newcommand{\amfgs}{\mbox{AMF-GS}}

\newlength\maxsteps

\definecolor{matlabblue}{HTML}{0072BD}
\definecolor{matlaborange}{HTML}{D95319}
\definecolor{matlabyellow}{HTML}{EDB120}
\definecolor{matlabpurple}{HTML}{7E2F8E}
\definecolor{matlabgreen}{HTML}{77AC30}
\definecolor{matlablightblue}{HTML}{4DBEEE}
\definecolor{matlabred}{HTML}{A2142F}

\tikzstyle{lines} = [
  solid,
  line width = 2pt
]

\tikzstyle{symbols} = [
  mark size    = 1.6pt,
  mark options = {solid}
]

\tikzstyle{gs} = [
  matlabblue,
  lines,
  dashed
]

\tikzstyle{rmfgs} = [
  matlaborange,
  lines,
  solid
]

\tikzstyle{amfgs} = [
  matlabpurple,
  lines,
  solid
]

\pgfplotscreateplotcyclelist{railplotlist}{
  {gs, symbols, mark = *},
  {rmfgs, symbols, mark = square*},
  {rmfgs, symbols, mark = diamond*},
  {rmfgs, symbols, mark = triangle*},
  {rmfgs, symbols, mark = pentagon*},
  {rmfgs, symbols, mark = *},
  {amfgs, symbols, mark = square*},
  {amfgs, symbols, mark = diamond*},
  {amfgs, symbols, mark = triangle*},
  {amfgs, symbols, mark = pentagon*},
  {amfgs, symbols, mark = *}
}

\pgfplotscreateplotcyclelist{cylinderplotlist}{
  {gs, symbols, mark = triangle*},
  {rmfgs, symbols, mark = square*},
  {rmfgs, symbols, mark = diamond*},
  {rmfgs, symbols, mark = triangle*},
  {amfgs, symbols, mark = square*},
  {amfgs, symbols, mark = diamond*},
  {amfgs, symbols, mark = triangle*}
}

\begin{document}
  

\title{Multi-fidelity robust controller design with gradient sampling}
  
\author[$\ast$,1]{Steffen W. R. Werner}
\affil[$\ast$]{Courant Institute of Mathematical Sciences,
  New York University, New York, NY 10012, USA.}
\affil[1]{\email{steffen.werner@nyu.edu}, \orcid{0000-0003-1667-4862}}
  
\author[$\ast$,2]{Michael L. Overton}
\affil[2]{\email{overton@cims.nyu.edu}, \orcid{0000-0002-6563-6371}}

\author[$\ast$,3]{Benjamin Peherstorfer}
\affil[3]{\email{pehersto@cims.nyu.edu}, \orcid{0000-0002-1558-6775}}
  
\shorttitle{Multi-fidelity gradient sampling}
\shortauthor{S.~W.~R. Werner, M.~L. Overton, B. Peherstorfer}
\shortdate{2022-12-05}
\shortinstitute{}
  
\keywords{%
  nonsmooth optimization,
  multi-fidelity methods,
  robust control,
  linear dynamical systems,
  H-infinity norm
}

\msc{%
  37N35, 
  37N40, 
  65K10, 
  90C30, 
  90C59  
}
  
\abstract{%
  Robust controllers that stabilize dynamical systems even under disturbances
  and noise are often formulated as solutions of nonsmooth, nonconvex
  optimization problems.
  While methods such as gradient sampling can handle the nonconvexity and
  nonsmoothness, the costs of evaluating the objective function may be
  substantial, making robust control challenging for dynamical systems with
  high-dimensional state spaces.
  In this work, we introduce multi-fidelity variants of gradient sampling that
  leverage low-cost, low-fidelity models with low-dimensional state spaces for
  speeding up the optimization process while nonetheless providing convergence
  guarantees for a high-fidelity model of the system of interest, which is
  primarily accessed in the last phase of the optimization process.
  Our first multi-fidelity method initiates gradient sampling on higher
  fidelity models with starting points obtained from cheaper, lower
  fidelity models.
  Our second multi-fidelity method relies on ensembles of gradients that
  are computed from low- and high-fidelity models.
  Numerical experiments with controlling the cooling of a steel rail profile and
  laminar flow in a cylinder wake demonstrate that our new multi-fidelity
  gradient sampling methods achieve up to two orders of magnitude speedup
  compared to the single-fidelity gradient sampling method that relies on the
  high-fidelity model alone.
}

\novelty{}

\maketitle


\section{Introduction}%
\label{sec:intro}

Robust controllers are a ubiquitous tool to overcome uncertainties in the
control of real-world applications resulting from the gap between mathematical
modeling and reality.
Constructing such controllers via minimizing the $\Hinf$-norm of closed-loop
systems is numerically challenging for at least two reasons.
First, the optimization objective induced by $\Hinf$-control leads to a
challenging optimization problem due to its nonsmooth and nonconvex nature.
Second, each evaluation of the objective entails computing the $\Hinf$-norm,
which incurs costs that grow rapidly with the dimension of the state space of
the system model. Gradient sampling methods~\cite{BurCLetal20, BurLO05, Kiw07}
can handle the nonsmooth, nonconvex objectives underlying $\Hinf$-control;
however, each evaluation of the objective remains computationally expensive.
We introduce multi-fidelity approaches that build on gradient sampling and
leverage hierarchies of low-fidelity models of the system of interest for
speeding up the optimization while still providing
convergence guarantees for the high-fidelity model
of the system.
Our new multi-fidelity variants of gradient sampling make finding
$\Hinf$-controllers tractable for models of systems with high-dimensional state
spaces, where relying on the expensive high-fidelity model alone quickly becomes
computationally prohibitive.
  
Multi-fidelity methods for optimization have a long tradition, especially in the
engineering community; see, e.g., the survey~\cite{PehWG18}.
Early work on multi-fidelity optimization 
was based on trust-region
methods~\cite{AleDLetal98, AriFS02, FahS02, FisGB17, RobEWetal12}.
Other works use a combination of reduced and full models in
optimization~\cite{PehW15a, PehW15, QiaGVetal17, ZahF14} and especially target
optimization under uncertainty, where the objective depends on stochastic
auxiliary variables~\cite{ChaPW20, KraPW17, LawCP22, NgW14, NgW16, PehWG16}.
For optimization problems with constraints given by partial differential
equations (PDEs), e.g., optimal control problems with smooth objective
functions, hierarchies of discretizations of PDEs have been used for efficient
preconditioning~\cite{BorK05, FahS02, HerS10, PeaSW14}.
In the context of uncertainty quantification, warm-starting iterative processes
is a common multi-fidelity approach; see, e.g.,~\cite{AlsVP22}.
There are also derivative-free multi-fidelity
methods~\cite{MarW12, MarW12a, WilS11}; however, these still require a smooth
objective function and thus are not well-suited for nonsmooth optimization
problems arising in $\Hinf$-control.

There is a large body of work on reduced modeling for control and control for
large-scale systems; see, e.g.,~\cite{Ben04, BenLP08, ObiA01, ReiS08a}.
The problem of efficiently designing $\Hinf$-controllers for
large-scale systems has been addressed before from different view points.
While in~\cite{MitO15} a new large-scale $\Hinf$-norm computation routine
was used to improve performance of optimization algorithms, reduced-order
surrogates were instead exploited in~\cite{BenMO18}.
In~\cite{BenHW22,MusG91}, analytical formulas for (suboptimal)
$\Hinf$-controllers are used rather than an optimization algorithm, relating the
low-order controller design problem under additional assumptions to the solution
of large-scale sparse nonlinear matrix equations.
  
The multi-fidelity variants of gradient sampling that we introduce in this work
can cope with nonconvex, nonsmooth objectives and at the same time leverage
low-fidelity models for reducing the optimization costs.
In the first multi-fidelity method that we introduce, we start by
optimizing the objective corresponding to a low-fidelity model and then use the
last iterate from the lower level as a starting point for optimization of the
objective corresponding to the next level.
This process is repeated until we eventually optimize with respect to the most
expensive, high-fidelity model with a good starting point.
The second variant uses the high-fidelity model to compute the objective
function and its gradient throughout the calculation, but restricts the
typically expensive gradient sampling process to gradients of the
lower-fidelity models until the final phase of the computation.
Numerical experiments demonstrate that speedups of up to two orders of magnitude
can be obtained compared to single-fidelity gradient sampling that uses the
high-fidelity model alone.
  
The paper is organized as follows.
We first discuss $\Hinf$-control and gradient sampling
methods in \Cref{sec:basics}.
We then introduce two new multi-fidelity variants of gradient sampling in
\Cref{sec:methods}.
We present numerical experiments for both variants on two real-world
applications, controlling the cooling of a steel rail profile and control of a
laminar flow in a cylinder wake in \Cref{sec:examples}.
Conclusions are drawn in \Cref{sec:conclusions}.


\section{Mathematical preliminaries}%
\label{sec:basics}

This section reviews the concepts of linear state-space systems, robust
$\Hinf$-controller design and the gradient sampling method.


\subsection{Dynamical systems and feedback controllers}

Consider a finite-di\-men\-sion\-al open-loop state-space model
of the form
\begin{equation} \label{eqn:sys}
  G\colon \begin{cases}
    \begin{aligned}
      E \dot{x}(t) & = _{\phantom{1}}A x(t) +
        _{\phantom{1}}B_{1} w(t) + _{\phantom{1}}B_{2} u(t),\\
      z(t) & = C_{1} x(t) + D_{11} w(t) + D_{12} u(t),\\
      y(t) & = C_{2} x(t) + D_{21} w(t) + D_{22} u(t),
    \end{aligned}
  \end{cases}
\end{equation}
where $x(t) \in \R^{\nG}$ are the internal states,
$u(t) \in \R^{\mG_{2}}$ the control inputs,
$w(t) \in \R^{\mG_{1}}$ the disturbances,
$z(t) \in \R^{\pG_{1}}$ the performance of the system and
$y(t) \in \R^{\pG_{2}}$ the measurements.
The matrices describing the model have corresponding dimensions:
$E, A \in \R^{\nG \times \nG}$,
$B_{1} \in \R^{\nG \times \mG_{1}}$,
$B_{2} \in \R^{\nG \times \mG_{2}}$,
$C_{1} \in \R^{\pG_{1} \times \nG}$,
$C_{2} \in \R^{\pG_{2} \times \nG}$,
$D_{11} \in \R^{\pG_{1} \times \mG_{1}}$,
$D_{12} \in \R^{\pG_{1} \times \mG_{2}}$,
$D_{21} \in \R^{\pG_{2} \times \mG_{1}}$ and
$D_{22} \in \R^{\pG_{2} \times \mG_{2}}$;
see, e.g.,~\cite{Fra87,ZhoD99}.
The system structure of~\cref{eqn:sys} is motivated by the observation that
mathematical models are inevitably idealized and that allowance must be made for
perturbations to the system, either because of its complexity in practice or
because of unpredictable external input.
The system~\cref{eqn:sys} therefore has two different types of inputs:
a deterministic signal $u$ that is the output of a controller, and a second
signal $w$ that accounts for modeling errors and random perturbations.
Furthermore, \cref{eqn:sys} has two outputs, one called $y$ that represents
state measurements, typically obtained by sensors, and a second output $z$,
which may not be measured in practice but represents the overall performance of
the system.
We consider~\cref{eqn:sys} without any direct feed-through term,
i.e., $D_{22} = 0$, to simplify the exposition.
In the general case with $D_{22} \neq 0$, it is described
in~\cite[Sec.~14.7]{ZhoD99} how one may first construct a controller $K$ with
transfer function $K(s)$ for the system with $D_{22} = 0$ and 
then obtain the controller for the system with $D_{22} \neq 0$ from
$K(s) \big( I_{p_{2}} + D_{22} K(s) \big)^{-1}$.
Also, we assume the matrix pencil $\lambda E - A$ in~\cref{eqn:sys}
to be regular, i.e., there exists a $\lambda \in \C$ such that $\lambda E - A$
is invertible, so that~\cref{eqn:sys} has a classical frequency domain
representation in terms of a transfer function.

The goal is to construct a continuous-time, finite-dimensional, feedback
controller, which maps the measurements taken from~\cref{eqn:sys} onto an
appropriate control signal, $K\colon y \mapsto u$.
The controller takes the form of a linear state-space model with
\begin{equation} \label{eqn:K}
  K\colon \begin{cases}
    \begin{aligned}
      \dxK(t) & = \AK \xK(t) + \BK y(t),\\
      u(t) & = \CK \xK(t) + \DK y(t),
    \end{aligned}
  \end{cases}
\end{equation}
where $\AK \in \R^{\nK \times \nK}$, $\BK \in \R^{\nK \times \pG_{2}}$,
$\CK \in \R^{\mG_{2} \times \nK}$ and $\DK \in \R^{\mG_{2} \times \pG_{2}}$.
Here, $\nK \in \N$ is the order of the controller, assumed
to be a fixed number that is much smaller than the state-space dimension $\nG$
of the system to be controlled, so $\nK \ll \nG$.
Note that, in contrast to the open-loop system~\cref{eqn:sys},
the controller~\cref{eqn:K} does not have a descriptor
(mass) matrix $\EK$; this
is motivated by engineering practice that avoids the use of active
algebraic constraints in the controller.
The control loop of~\cref{eqn:sys} is  closed by connecting the
controller~\cref{eqn:K} with the system~\cref{eqn:sys}, which yields the
closed-loop system $\Gclosed\colon w \mapsto z$ with
\begin{equation} \label{eqn:closedsys}
  \Gclosed\colon \begin{cases}
    \begin{aligned}
      \Eclosed \dot{x}_{\mathsf{c}} & = \Aclosed x_{\mathsf{c}} +
        \Bclosed w(t),\\
      z(t) & = \Cclosed x_{\mathsf{c}} + \Dclosed w(t),
    \end{aligned}
  \end{cases}
\end{equation}
where the system matrices are given by
\begin{equation} \label{eqn:EcAcBcCcDc}
  \begin{aligned}
    \Eclosed & = \begin{bmatrix} E & 0 \\ 0 & I_{\nK} \end{bmatrix}, &
    \Aclosed & = \begin{bmatrix} A + B_{2} \DK C_{2} & B_{2} \CK \\
      \BK C_{2} & \AK \end{bmatrix}, \\
    \Bclosed & = \begin{bmatrix} B_{1} + B_{2} \DK D_{21} \\
      \BK D_{21} \end{bmatrix}, &
    \Cclosed & = \begin{bmatrix} C_{1} + D_{12} \DK C_{2} & D_{12}
      \CK \end{bmatrix}, \\
    \Dclosed & = D_{11} + D_{12} \DK D_{21}.
  \end{aligned}
\end{equation}


\subsection{\texorpdfstring{$\Hinf$}{H-infinity} controller design}%
\label{subsec:hinfcon}

The requirement for the feedback controller~\cref{eqn:K} that we consider
here is the stabilization of the closed-loop system~\cref{eqn:closedsys}, i.e.,
the design of~\cref{eqn:K} ensures that the closed-loop matrix pencil
$s \Eclosed - \Aclosed$ is regular and that all of its finite eigenvalues lie in
the open left half-plane.
Thus, we define the set of stabilizing controllers as
\begin{equation*}
  \K = \left\{ (\AK, \BK, \CK, \DK) \left|~
    \lambda \in \C ~\text{with}~ \det(\lambda \Eclosed - \Aclosed) = 0
    ~\Rightarrow~ \real(\lambda) < 0
    \right.\right\}.
\end{equation*}
Let $\lVert \cdot \rVert_{\Hinf}$ denote the $\Hinf$-norm, defined
for the closed-loop system~\cref{eqn:closedsys} by
\begin{equation*}
  \lVert \Gclosed \rVert_{\Hinf} := \sup\limits_{\lambda \in \C,
    \real(\lambda) \geq 0} \lVert \Gclosed(\lambda) \rVert_{2},
\end{equation*}
with the transfer function $\Gclosed(s) = \Cclosed (s\Eclosed -
\Aclosed)^{-1} \Bclosed + \Dclosed$, where $s \in \C$; see, e.g.,~\cite{Ant05}.
In optimal $\Hinf$-control, a controller $K \in \K$ is sought as a solution to
the constrained minimization problem
\begin{equation} \label{eqn:HinfK}
  \min\limits_{K \in \K} \lVert \Gclosed \rVert_{\Hinf}.
\end{equation}
The task of $\Hinf$-optimal control can be interpreted as finding a stabilizing
controller that minimizes the worst-case amplification of all admissible
disturbances.

In this paper, we focus on the case where the open-loop
system~\cref{eqn:sys} and, consequently, the closed-loop
system~\cref{eqn:closedsys}, are described by large-scale sparse systems of
differential-algebraic equations.
The spectral abscissa of the pencil $s \Eclosed - \Aclosed$ is the real part of
its rightmost finite eigenvalue; we denote this by
\begin{equation} \label{eqn:specabs}
  \alpha(\Aclosed, \Eclosed) := \max \left\{ \real(\lambda) \left|~
    \lambda \in \C ~\text{with}~
    \det(\lambda \Eclosed - \Aclosed) = 0 \right. \right\}.
\end{equation}
The maximum peak of the spectral norm of the transfer function on the
imaginary axis is known as the $\Linf$-norm, which is for the
closed-loop system~\cref{eqn:closedsys} given by
\begin{equation} \label{eqn:linfnorm}
  \lVert \Gclosed \rVert_{\Linf} := \sup\limits_{\omega \geq 0}
    \lVert \Gclosed(\i\,\omega) \rVert_{2},
\end{equation}
where $\i$ denotes the imaginary unit, and the supremum is over the nonnegative
imaginary axis because the data are real.

Using~\cref{eqn:specabs,eqn:linfnorm}, the $\Hinf$-norm is
\begin{equation} \label{eqn:hinfnorm2}
  \lVert \Gclosed \rVert_{\Hinf} =
    \begin{cases}
      \lVert \Gclosed \rVert_{\Linf} &
        \text{if}~\alpha(\Aclosed, \Eclosed) < 0,\\
      \infty & \text{otherwise}.
    \end{cases}
\end{equation}
Now we define our objective function to be minimized as
\begin{equation} \label{eqn:fdef}
       f(x) := \lVert \Gclosed \rVert_{\Hinf}
\end{equation}
with the design variable
\begin{equation} \label{eqn:designvar}
  x = \begin{bmatrix} \vec(\AK) \\ \vec(\BK) \\ \vec(\CK) \\ \vec(\DK)
    \end{bmatrix}\in \R^{\nopt} ~\text{where}~
    \nopt = \nK^{2} + \nK m_{2} + p_{2} \nK + p_{2} m_{2},
\end{equation}
defining a controller~\cref{eqn:K} via the matrices $K = (\AK, \BK, \CK, \DK)$,
with the closed-loop system matrices defining $\Gclosed$ in~\cref{eqn:hinfnorm2}
depending on $K$ via~\cref{eqn:EcAcBcCcDc}.
It is also convenient to define the constraint function
\begin{equation} \label{eqn:hdef}
  h(x) := \alpha(\Aclosed, \Eclosed),
\end{equation}
where again the closed-loop system matrices $\Aclosed$ and $\Eclosed$ depend 
on the controller matrices $K = (\AK, \BK, \CK, \DK)$ via~\cref{eqn:EcAcBcCcDc}.
Using this notation, the optimization problem~\cref{eqn:HinfK} may be
equivalently given as either
\begin{equation} \label{eqn:HinfX}
  \begin{aligned}
    \min\limits_{x} f(x) && \text{or} &&
    \min\limits_{x:h(x)<0} f(x).
  \end{aligned}
\end{equation}
This optimization problem is challenging because the
$\Hinf$-norm~\cref{eqn:hinfnorm2} is nonconvex and, at points $x$ where the
supremum in~\cref{eqn:linfnorm} is attained at more than one value of $\omega$,
nonsmooth.
However, $f$ is locally Lipschitz on the set of stabilizing controllers
$\{x \in \R^{\nopt}: h(x) < 0 \}$.


\subsection{Gradient sampling method}%
\label{subsec:gradsamp}

It has been known for decades that the steepest descent method (gradient
descent with a line search) generally fails on nonsmooth optimization problems,
typically converging to a non-stationary (and non-optimal) point where the
objective function is not differentiable. 
The gradient sampling method is a stabilized steepest descent method devised to
overcome this difficulty. It was presented by Burke, Lewis and Overton
in 2005~\cite{BurLO05}, along with an extensive convergence theory that was
subsequently refined by Kiwiel in 2007~\cite{Kiw07}. 
The algorithm is nondeterministic in the sense that it generates (samples)
gradients at randomly generated points within an appropriately-sized ball
around a given iterate.
In this paper, we rely on the detailed description of the method and its
convergence theory in the survey~\cite{BurCLetal20}.
The main convergence result for Alg.~GS of~\cite{BurCLetal20} (with specific
parameter choices) is stated as Theorem~6.1 there:
\emph{Suppose that $f$ is locally Lipschitz on $\R^{\nopt}$ and continuously
differentiable on an open set with full measure.
Then, with probability one, Alg.~GS is well defined and does not terminate,
and generates a sequence of iterates for which either the function values
diverge to $-\infty$, or every cluster point of the sequence is Clarke
stationary for $f$.}
Clarke stationarity is a standard measure of stationarity for locally Lipschitz,
nonsmooth functions~\cite{BorL06}.

The gradient sampling method relies on the computation of the function $f$ and
its gradient $\grad f$ at the sequence of iterates generated by the method,
using a ``gradient paradigm''~\cite{AslO21}, as opposed to the``subgradient 
paradigm'' often used for nonsmooth functions, in particular by the
``subgradient method'', which is usually very slow.
The gradient paradigm observes that, since locally Lipschitz functions are
differentiable almost everywhere by Rademacher's theorem, and since in practice,
it is essentially impossible to verify whether a nontrivial function $f$ is
differentiable or not at a given iterate $x$, a method can reasonably compute an
approximate gradient at any given point, for example, by ignoring ``ties'' in a
max function.
The idea is that it is only in the limit of the sequence of iterates that the
function is actually not differentiable.
Of course, sampled gradients computed at nearby points in this way may vary
greatly, and the gradient sampling algorithm exploits this property.
These key points are discussed at greater length in the references given above.

The gradient sampling method has been applied to solve $\Hinf$-norm optimization
and related stabilization problems since it was first
introduced~\cite{BurHLetal06,BurHLetal06a,BurLO03}.
We follow the same basic strategy used in~\cite{BurHLetal06}:
First, in order to find a stabilizing controller for the $\Hinf$-norm
optimization problem described in \Cref{subsec:hinfcon}, we apply gradient
sampling to the constraint function $h(x)$ defined in~\cref{eqn:hdef};
then, once a point $x^{0}$ with $h(x^{0})<0$ has been found, we apply gradient
sampling to the $\Hinf$-norm objective $f$ defined in~\cref{eqn:fdef},
initialized at $x^{0}$.
If this results in $f$ being evaluated at a non-stabilizing controller, the
function value $\infty$ that is returned will result in the controller being
rejected by the line search; according to the gradient sampling convergence
theory, as long as $f$ is differentiable at $x^{0}$, the line search must
eventually return a new point $x^{1}$ with $f(x^{1}) < f(x^{0})$.  
The functions $f$ and $h$ are differentiable almost everywhere
(in the former case, almost everywhere on $\mathcal{K}$), and the
formulas for their gradients may be derived from the formulas for
the gradients of the $\Hinf$-norm and the spectral abscissa given in
\Cref{app:hinfnorm,app:specabs}, respectively.


\section{Multi-fidelity gradient sampling}%
\label{sec:methods}

In this section, we introduce two multi-fidelity versions of the gradient
sampling method to design controllers for high-fidelity models for which
a hierarchy of cheap low-fidelity models are available.
We first introduce the notation of hierarchies of models in
\Cref{sec:ML:HierModels}.
Then we define our two new methods:
Gradient sampling with multi-fidelity restarts in \Cref{subsec:rmfgs} and
gradient sampling with multi-fidelity approximate gradients
in \Cref{subsec:amfgs}.


\subsection{Hierarchies of models}
\label{sec:ML:HierModels}

We consider the situation where there is a hierarchy of $\maxdisc$ models
of the form~\cref{eqn:sys} available.
The accuracy of the models increases with a corresponding index from level~$1$
to level~$\maxdisc$, the most accurate model.
We find such a situation, for example, when~\cref{eqn:sys} is given as
spatial discretization of partial differential equations, where the
model hierarchy with levels $\ell = 1, \dots, \maxdisc$ is due to different
refinements of the discretization.
The hierarchy of models gives rise to a hierarchy of objective functions for
$\Hinf$-controller design:
\begin{equation} \label{eqn:hinflevel}
  f^{\ell}(x) = \lVert \Gclosed^{\ell} \rVert_{\Hinf},
\end{equation}
with $\ell = 1, \ldots, \maxdisc$.
A key point to note is that the dimension $\nopt$ of the vector $x$
in~\cref{eqn:hinflevel} representing the  controller $K = (\AK, \BK, \CK, \DK)$
is independent of the model level~$\ell$.
Instead of~\cref{eqn:EcAcBcCcDc}, we now have closed-loop system matrices
defined by
\begin{equation*}
  \begin{aligned}
    \Eclosed^{\ell} & = \begin{bmatrix} E^{\ell} & 0 \\ 0 & I_{\nK}
      \end{bmatrix}, &
    \Aclosed^{\ell} & = \begin{bmatrix} A^{\ell} + B_{2}^{\ell} \DK
      C_{2}^{\ell} & B_{2}^{\ell} \CK \\
      \BK C_{2}^{\ell} & \AK \end{bmatrix}, &
    \Bclosed^{\ell} & = \begin{bmatrix} B_{1}^{\ell} + B_{2}^{\ell} \DK
      D_{21}^{\ell} \\
      \BK D_{21}^{\ell} \end{bmatrix},\\
    & &
    \Cclosed^{\ell} & = \begin{bmatrix} C_{1}^{\ell} + D_{12}^{\ell} \DK
      C_{2}^{\ell} & D_{12}^{\ell} \CK \end{bmatrix}, &
    \Dclosed^{\ell} & = D_{11}^{\ell} + D_{12}^{\ell} \DK D_{21}^{\ell},
  \end{aligned}
\end{equation*}
where the matrices superscripted by $\ell$ are the open-loop system matrices.
The corresponding transfer functions of the closed-loop systems are
$\Gclosed^{\ell}(s) = \Cclosed^{\ell} (s\Eclosed^{\ell} -
\Aclosed^{\ell})^{-1} \Bclosed^{\ell} + \Dclosed^{\ell}$.

Our aim is to find a controller that is optimal with respect to the
high-fidelity objective function $f^{\maxdisc}$, while leveraging the less
accurate but cheaper objective functions $f^{\ell}$ on levels
$\ell = 1, \ldots, \maxdisc - 1$.
The objective functions have gradients $\grad f^{1}, \dots, \grad f^{\maxdisc}$,
which are increasingly more expensive to compute as $\ell$ increases;
see \Cref{app:hinfnorm} for the formulas.

Besides hierarchies of discretizations, the model hierarchy may alternatively be
obtained via model reduction techniques. 
These allow the computation of reasonably accurate, cheap-to-evaluate
surrogates that can serve as low-fidelity models in our setting.
See, for example,~\cite{BenCOetal17, BenGW15, BenSGetal21, BenSGetal21a, QuaR14}
for overviews on potential methods, or~\cite{BenHW22, MusG91} for model
reduction methods in the context of $\Hinf$-controller design.

In the following, our starting point is a hierarchy of objective functions
$f^{1}, \dots, f^{\maxdisc}$ that are ordered from cheap to expensive and less
accurate to more accurate but we make no assumptions on where the objective
functions originate.
It is sufficient to have an oracle that allows the evaluation of the functions
$f^{\ell}$ and their gradients $\grad f^{\ell}$ at the design variable $x$
corresponding to the given controller $K$.
Besides hierarchies of objective functions, we must also at least implicitly
consider hierarchies of constraint functions $h^{\ell}(x)$.
We return to this topic below.


\subsection{Restarted multi-fidelity gradient sampling (\rmfgs{})}%
\label{subsec:rmfgs}

Our restart\-ed multi-fidelity gradient sampling (\rmfgs{}) method uses
controllers obtained with lower fidelity models to warm-start the
optimization for controllers of higher fidelity models.


\subsubsection{Multi-fidelity restarts}

The proposed \rmfgs{} approach iterates over the levels
$\ell = 1, \dots, \maxdisc$ and, at each level $\ell$, solves an
optimization problem of the form~\cref{eqn:HinfX} with the objective function
$f^{\ell}$, where the initial guess is the solution of the previous level. So,
letting $x^{k_{\ell-1}}$ denote the final iterate at level $\ell - 1$, 
the initial guess at level $\ell \geq 2$ is $x^{k_{\ell-1}}$. The motivation for
\rmfgs{} is that the objective functions become progressively
more accurate with increasing level $\ell$, and thus, the solution
$x^{k_{\ell-1}}$ at the previous level $\ell - 1$ should be a good
starting point at the current level $\ell$, implying that fewer gradient
sampling steps are necessary than with a generic initial guess.
Hence, the aim is to take many iterations on lower levels where the initial
starting points are poor but where objective and gradient evaluations are cheap,
while taking fewer of the expensive evaluations on higher levels as the starting
points get closer to a minimizer of the high-fidelity objective
function~$f^{\maxdisc}$.

For any level $\ell$, the function $f^{\ell}$ is monotonically decreasing on
$\{ x^{k} \}$ as $k$ increases from $k_{\ell-1}$ to $k_\ell$.
Note, however, that for $\ell < \maxdisc$, there is no guarantee that the
high-fidelity objective $f^{\maxdisc}$ is lower at
$x^{k_\ell}$ than it was at $x^{k_{\ell-1}}$.
Indeed, it might not even be finite, since the objective function is finite only
if the closed-loop system is stable, and even if this is the case for the model
at one level, it might not be at another level.


\subsubsection{Algorithmic description of \rmfgs{}}

\begin{algorithm}[t]
  \caption{Restarted multi-fidelity gradient sampling (\rmfgs).}
  \label{alg:rmfgs}
  \hspace*{0\baselineskip}\textbf{Input:}~Initial point $x^{0} \in \R^{\nopt}$,
    \par\phantom{\textbf{Input:}}
    sample size 
    $ q \geq \nopt + 1$,
    initial sampling radii $\epsilon_{\ell, 0} > 0$,
    \par\phantom{\textbf{Input:}}
    initial stationarity targets $\nu_{\ell, 0} > 0$,
    \par\phantom{\textbf{Input:}}
    termination tolerances $\epsilon_{\ell, \opt} \in (0, \epsilon_{\ell, 0}), 
    \nu_{\ell, \opt} \in (0, \nu_{\ell, 0})$,
    \par\phantom{\textbf{Input:}}
    reduction factors $\theta_{\ell, \epsilon} \in (0,1), 
    \theta_{\ell, \nu} \in (0, 1)$, and
    \par\phantom{\textbf{Input:}}
    line search parameters $\beta_{\ell} \in (0,1)$, $\gamma_{\ell} \in (0,1)$,
    for $\ell = 1, \ldots, \maxdisc$.
      
  \hspace*{0\baselineskip}\textbf{Output:}~Approximation $x^{k} \in \R^{\nopt}$
    to a minimizer of $f^{\maxdisc}$.

  \begin{algorithmic}[1]
    \State Initialize $k = 0$.
    \For{$\ell = 1$ \textbf{to} $\maxdisc$} \label{alg:rmfgs:outloopstart}
      \If{$f^{\ell}(x^{k})$ is not finite}
        \State Apply stabilization step for $f^{\ell}$ to $x^{k}$.
      \EndIf
      \State Set $\nu_{k+1} = \nu_{\ell, 0}$ and
        $\epsilon_{k+1} = \epsilon_{\ell, 0}$.
      \Repeat \label{alg:rmfgs:inloopstart}
        \State Independently sample $\{ x^{k, 1}, \ldots, x^{k, q}\}$ uniformly
          from $\mathcal{B}(x^{k}, \epsilon_{k})$.
          \label{alg:rmfgs:sample}
        \State Compute $g^{k}$ as the solution of $\min_{g \in
          \mathcal{G}^{\ell, k}} \frac{1}{2} \lVert g \rVert_{2}^{2}$, where
          \begin{equation*}
            \mathcal{G}^{\ell, k} = \conv\left\{ \nabla f^{\ell}(x^{k}), 
              \nabla f^{\ell}(x^{k, 1}), \ldots, \nabla f^{\ell}(x^{k, q})
              \right\}.
          \end{equation*}
          \label{alg:rmfgs:grad}
          \vspace{-\baselineskip}
        \State Compute $x^{k + 1}, \epsilon_{k + 1}, \nu_{k + 1}$ using
          \Cref{alg:gsstep} with inputs
          \par\hspace{.5\baselineskip}
          $x^{k}, g^{k}, f^{\ell}, \epsilon_{k}, \nu_{k},
          \theta_{\ell, \epsilon}, \theta_{\ell, \nu},
          \epsilon_{\ell, \opt}, \nu_{\ell, \opt}, \beta_{\ell}, \gamma_{\ell}$.
          \label{alg:rmfgs:step}
        \State Increment $k \leftarrow k + 1$.
      \Until{($x^{k} == x^{k-1}$) \textbf{and}
        ($\epsilon_{k} == \epsilon_{k-1}$) \textbf{and}
        ($\nu_{k} == \nu_{k-1}$).} \label{alg:rmfgs:inloopend}
    \EndFor \label{alg:rmfgs:outloopend}
  \end{algorithmic}
\end{algorithm}

The new method is summarized in \Cref{alg:rmfgs}.
The main difference from the original (single-fidelity) gradient sampling
method~\cite[Alg.~GS]{BurCLetal20} is the new outer loop starting in
Line~\ref{alg:rmfgs:outloopstart} of \Cref{alg:rmfgs}, which iterates over the
available levels $\ell = 1, \ldots, \maxdisc$.
Lines~\ref{alg:rmfgs:inloopstart} to~\ref{alg:rmfgs:inloopend} consist of an
inner iteration describing the single-fidelity
gradient sampling method using the objective function~$f^{\ell}$ and
its gradients~$\grad f^{\ell}$ at the current level.
This has three parts:
\begin{itemize}
  \item[(a)] In Line~\ref{alg:rmfgs:sample}, sampling gradients uniformly from
    $\mathcal{B}(x^{k}, \epsilon_{k})$, the 2-norm ball around the current
    iterate $x^{k}$ with radius~$\epsilon_{k}$.
  \item[(b)] In Line~\ref{alg:rmfgs:grad}, computing the vector $g^{k}$, which
    is easily done by standard software for convex quadratic programming,
    observing that the convex hull of vectors
    $v^{1}, \ldots, v^{q} \in \R^{\nopt}$ is 
    \begin{equation*}
      \left\{ \alpha_{1} v^{1} + \ldots + \alpha_{q} v^{q} ~\left|~
        \alpha_{1} + \ldots + \alpha_{q} = 1,
        \alpha_{1} \geq 0, \ldots, \alpha_{q} \geq 0 \right. \right\}.
    \end{equation*}%
    As explained in~\cite[Sec.~6.1]{BurCLetal20},
    the vector $-g^{k}$ is not only a descent direction for
    $f^{\ell}$, but more importantly it is a \emph{stabilized} or \emph{robust}
    descent direction, which allows for longer steps to be taken in the line
    search in the next part.
  \item[(c)] In Line~\ref{alg:rmfgs:step}, the computation of the gradient
    sampling step as
    described in \Cref{alg:gsstep}, which includes checking the convergence
    criteria, updating the algorithm parameters accordingly, and, if the
    termination criteria are not yet met, updating the current iterate using a
    line search along $-g^{k}$.
\end{itemize}

The inner iteration for a given $f^{\ell}$ terminates when the gradient sampling
step has no effect, i.e., if the new iterate is the same as the previous one
and the sampling radius and stationarity target did not change.
Looking at \Cref{alg:gsstep}, we see that this can only occur if the algorithm
satisfies the convergence criteria specified by the parameters.
According to the gradient sampling theory, this must happen eventually;
see~\cite[Cor.~6.1]{BurCLetal20}, taking into account
the initialization of the parameters in \Cref{alg:rmfgs}.
In practice, it is necessary to set a limit on the number of steps in each inner
iteration, both because of the possible effects of rounding errors and to limit
the overall computation time.
Likewise, in theory, the line search in Line~\ref{alg:gsstep:linesearch} of
\Cref{alg:gsstep} must terminate in a finite number of steps, although in
practice, because of rounding errors, a limit must be placed on this and the
line search terminated if this limit is reached.
Whichever way the iteration for level $\ell < \maxdisc$ terminates, 
the method continues with the next model level in the outer loop.
In this case, the current iterate $x^{k}$ is the final iterate $x^{k_{\ell}}$ of
level $\ell < \maxdisc$ and the initial iterate of level $\ell + 1$.

The algorithm allows for its parameters to depend on the
level $\ell$ so that adjustments for each level are possible.
The last step of the outer loop in \Cref{alg:rmfgs} is gradient sampling with
the objective function of interest $f^{\maxdisc}$, i.e., each step of the
inner loop in \Cref{alg:rmfgs} is as expensive as each step of classical
single-fidelity gradient sampling.
In terms of global computational costs in comparison to the single-fidelity
method~\cite[Alg.~GS]{BurCLetal20}, we can potentially save
function as well as gradient evaluations using \Cref{alg:rmfgs}, under the
assumption that the computed approximations of minimizers on each level are
indeed good initial guesses for optimization on subsequent levels.

\begin{algorithm}[t]
  \caption{Gradient sampling step.}
  \label{alg:gsstep}
  
  \hspace*{0\baselineskip}\textbf{Input:}~Iterate $x\in\R^{\nopt}$,
    vector $g \in\R^{\nopt}$,
    objective function $f$,
    \par\phantom{\textbf{Input:}}
    current sampling radius $\epsilon$ and stationarity target $\nu$,
    \par\phantom{\textbf{Input:}}
    reduction factors $\theta_{\epsilon}$ and $\theta_{\nu}$,
    termination tolerances $\epsilon_{\opt}$ and $ \nu_{\opt}$, and
    \par\phantom{\textbf{Input:}}
    line search parameters $\beta$ and $\gamma$.
      
  \hspace*{0\baselineskip}\textbf{Output:}~Updated iterate $\hat{x}$,
    sampling radius $\hat{\epsilon}$ and stationarity target $\hat{\nu}$.
    
  \begin{algorithmic}[1]
    \If{($\lVert g \rVert_{2} \leq \nu_{\opt}$) \textbf{and}
      ($\epsilon \leq \epsilon_{\opt}$)}
      \State Set $\hat{\nu} = \nu$,
        $\hat{\epsilon} = \epsilon$ and
        $\hat{t} = 0$.
    \Else
      \If{$\lVert g \rVert_{2} \leq \nu$}
        \State Set $\hat{\nu} = \theta_{\nu} \nu$,
          $\hat{\epsilon} = \theta_{\epsilon} \epsilon$ and 
          $\hat{t} = 0$.
      \Else
        \State Set $\hat{\nu} = \nu$ and
          $\hat{\epsilon} = \epsilon$.
        \State Set $\hat t = \max \left\{ t \in \{1, \gamma, \gamma^{2},
          \ldots \}: 
          f(x - tg) < f(x)- \beta t \|g\|_{2}^{2} \right\}$.
          \label{alg:gsstep:linesearch}
      \EndIf
    \EndIf
    \State Update $\hat{x} = x - \hat t g$.
  \end{algorithmic}
\end{algorithm}

\Cref{alg:gsstep} implements the update step of gradient sampling and is
the same as in Alg.~GS in~\cite{BurCLetal20}, except for the differentiability
check of the objective function $f$ at the next iterate $\hat{x}$.
This check is needed in theory in order to be able to rigorously state the
convergence results in~\cite{BurCLetal20}, but in practice, with the inevitable
rounding errors incurred in floating point arithmetic, it makes little or no
sense to attempt it.
As already noted, our objective functions are differentiable almost everywhere,
and while encountering a point where the function is actually not differentiable
is not technically a probability zero event, it may be considered extremely
unlikely in practice.
This issue is discussed further in~\cite[Sec.~6.4.2]{BurCLetal20}.


\subsection{Approximate multi-fidelity gradient sampling (\amfgs{})}%
\label{subsec:amfgs}

A valid criticism of \Cref{alg:rmfgs} is that although our primary interest is
in minimizing the highest fidelity model $f^{\maxdisc}$, this does not enter the
computation until the gradient sampling algorithm has been run on all lower
fidelity objectives $f^{1}, f^{2}, \ldots, f^{\maxdisc-1}$.
Although we justified this by arguing that the final iterate for one level
should be a good starting point for the next level, an alternative viewpoint is
that we might want to involve the highest fidelity model $f^{\maxdisc}$ at
earlier stages of the computation.
This can be done efficiently by using $f^{\maxdisc}$ as the objective function
from the beginning, but replacing the expensive gradient sampling of
$f^{\maxdisc}$ by gradient sampling of the cheaper models
$f^{1}, f^{2}, \ldots, f^{\maxdisc - 1}$.


\subsubsection{Multi-fidelity ensembles of gradients}

In the \amfgs{} method, we retain the idea of an outer loop over all $\maxdisc$
levels, but, unlike in the \rmfgs{} method, we involve the high-fidelity
function $f^{\maxdisc}$ at every stage of the outer loop.
For this reason, we enforce the property that the high-fidelity function
$f^{\maxdisc}$ is monotonically decreasing on $\{ x^{k} \}$ as $k$ increases.
However, although we evaluate $f^{\maxdisc}$ at every iterate $x^{k}$, and in
the line search that produces these iterates, it is only at the final level
$\maxdisc$ that we actually sample $q \geq \nopt + 1$ gradients of the
high-fidelity function $f^{\maxdisc}$.
At all earlier levels, we sample gradients of lower fidelity functions instead.
Thus, we replace the definition
\begin{equation*}
  \mathcal{G}^{\ell, k} = \conv\left\{
    \nabla f^{\ell}(x^{k}), 
    \nabla f^{\ell}(x^{k, 1}), \ldots,
    \nabla f^{\ell}(x^{k, q})
    \right\}
\end{equation*}
in Line~\ref{alg:rmfgs:grad} of \Cref{alg:rmfgs} by
\begin{equation*}
  \mathcal{G}^{\ell, k} = \conv\left\{
    \grad f^{\maxdisc}(x^{k}),
    \grad f^{\ell}(x^{k, 1}), \ldots,
    \grad f^{\ell}(x^{k, q}) \right\}.
\end{equation*}


\subsubsection{Algorithmic description of \amfgs{}}

\begin{algorithm}[t]
  \caption{Approximate multi-fidelity gradient sampling (\amfgs).}
  \label{alg:amfgs}
  \hspace*{0\baselineskip}\textbf{Input:}~Initial point $x^{0} \in \R^{\nopt}$,
    \par\phantom{\textbf{Input:}}
    sample size 
    $q \geq \nopt + 1$,
    initial sampling radii $\epsilon_{\ell, 0} > 0$,
    \par\phantom{\textbf{Input:}}
    initial stationarity targets $\nu_{\ell, 0} > 0$,
    \par\phantom{\textbf{Input:}}
    termination tolerances $\epsilon_{\ell, \opt} \in (0, \epsilon_{\ell, 0}), 
    \nu_{\ell, \opt} \in (0, \nu_{\ell, 0})$, 
    \par\phantom{\textbf{Input:}}
    reduction factors $\theta_{\ell, \epsilon} \in (0,1), 
    \theta_{\ell, \nu} \in (0, 1)$, and
    \par\phantom{\textbf{Input:}}
    line search parameters $\beta_{\ell} \in (0,1)$, $\gamma_{\ell} \in (0,1)$,
    for $\ell = 1, \ldots, \maxdisc$.
      
  \hspace*{0\baselineskip}\textbf{Output:}~Approximation $x^{k} \in \R^{\nopt}$
    to a minimizer of $f^{\maxdisc}$.

  \begin{algorithmic}[1]
    \State Initialize $k = 0$.
    \If{$f^{\maxdisc}(x^0)$ is not finite}
        \State Apply stabilization step for $f^{\maxdisc}$ to $x^{0}$.
     \EndIf
    \For{$\ell = 1$ \textbf{to} $\maxdisc$}
      \State Set $\nu_{k+1} = \nu_{\ell, 0}$ and
        $\epsilon_{k+1} = \epsilon_{\ell, 0}$.
      \Repeat
        \State Independently sample $\{ x^{k, 1}, \ldots, x^{k, q}\}$ uniformly
          from $\mathcal{B}(x^{k}, \epsilon_{k})$.
        \State Compute $g^{k}$ as the solution of $\min_{g \in
          \mathcal{G}^{\ell, k}} \frac{1}{2} \lVert g \rVert_{2}^{2}$, where
          \begin{equation*}
            \mathcal{G}^{\ell, k} = \conv\left\{ \nabla f^{\maxdisc}(x^{k}), 
              \nabla f^{\ell}(x^{k, 1}), \ldots, \nabla f^{\ell}(x^{k, q})
              \right\}.
          \end{equation*}
          \label{alg:amfgs:grad}
          \vspace{-\baselineskip}
        \State Compute $x^{k + 1}, \epsilon_{k + 1}, \nu_{k + 1}$ using
          \Cref{alg:gsstep} with inputs
          \par\hspace{.5\baselineskip}
          $x^{k}, g^{k}, f^{\maxdisc}, \epsilon_{k}, \nu_{k},
          \theta_{\ell, \epsilon}, \theta_{\ell, \nu},
          \epsilon_{\ell, \opt}, \nu_{\ell, \opt}, \beta_{\ell}, \gamma_{\ell}$.
        \State Increment $k \leftarrow k + 1$.
      \Until{($x^{k} == x^{k-1}$) \textbf{and}
        ($\epsilon_{k} == \epsilon_{k-1}$) \textbf{and}
        ($\nu_{k} == \nu_{k-1}$).}
    \EndFor
  \end{algorithmic}
\end{algorithm}

The \amfgs{} method is summarized in \Cref{alg:amfgs}.
The basic structure of the algorithm is the same as that of \Cref{alg:rmfgs}.
However, a major difference between them is that in \amfgs{}, we are minimizing
the high-fidelity objective function $f^{\maxdisc}$ at \emph{all} levels
$\ell = 1, \ldots, \maxdisc$, while in \rmfgs{}, at level $\ell$,
we minimize the objective $f^{\ell}$.
Consequently, each step of level $\ell$ of \amfgs{}~(\Cref{alg:amfgs}) is
computationally more expensive than the corresponding step in
\rmfgs{}~(\Cref{alg:rmfgs}). However, for $\ell < \maxdisc$, it is less
expensive than a step at level $\maxdisc$ of either method due to the use of
cheaper-to-evaluate approximations in the gradient computations of the
sampled evaluation points in Line~\ref{alg:amfgs:grad} of \Cref{alg:amfgs}.
A key point, however, is that at the current iterate $x^{k}$, we use the
gradient of the high-fidelity objective function $f^{\maxdisc}$ in the
definition of $\mathcal{G}^{\ell, k}$, regardless of the level $\ell$ in the
outer loop.
This guarantees that $-g^{k}$ is a descent direction for $f^{\maxdisc}$,
although how ``robust'' of a descent direction it is depends on how well the
sampled gradients of $f^{\ell}$ approximate gradients of $f^{\maxdisc}$.
If the approximation is not very good, the result may be that the line search
needs to take a very short step to obtain a reduction in $f^{\maxdisc}$
along~$-g^{k}$. 
The main differences between \Cref{alg:rmfgs,alg:amfgs} are the definition of
$\mathcal{G}^{\ell, k}$ and that the
function we pass to \Cref{alg:gsstep} is $f^\ell$ in the first case and
$f^{\maxdisc}$ in the second case.
Note that both methods, \Cref{alg:rmfgs,alg:amfgs}, boil down to the
classical (single-fidelity) gradient sampling method
from~\cite[Alg.~GS]{BurCLetal20} in the last step of each outer loop, so
the rationale for both methods is ultimately to provide a good starting point
for this final optimization at level $\maxdisc$. 


\subsection{Stabilization}%
\label{subsec:stabilization}

As explained in \Cref{subsec:gradsamp}, in order to obtain initial points for
minimization of the $\Hinf$-norm objective, it may be necessary to first apply
gradient sampling to the stabilization constraint function.
Thus, in \Cref{alg:rmfgs}, in order to initiate gradient sampling optimization
of $f^{\ell}$ at step $\ell$ of the outer loop, it may be necessary to first
apply gradient sampling to the corresponding constraint function $h^{\ell}$.
This applies not only at level $1$, but at higher levels as well, because
there is no guarantee that at level $\ell > 1$, the function $f^{\ell}$ is
finite at the starting point $x^{k}$, even though $f^{\ell-1}$ is necessarily
finite there.
However, we note that this stabilization step at level $\ell > 1$ was never
needed in our computational results presented in \Cref{sec:examples}.
In contrast, for~\cref{alg:amfgs}, at most one initial stabilization is
necessary, to obtain a point $x^{0}$ where $f^{\maxdisc}$ is finite.


\subsection{Theoretical guarantees}%
\label{subsec:theory}

Provided step $\ell$ in the outer loop of \Cref{alg:rmfgs} is initiated at a
point where $f^{\ell}$ is finite and differentiable, and that $f^{\ell}$ is
also differentiable at subsequent iterates (see the discussion
at the end of \Cref{subsec:rmfgs}), the convergence theory given
in~\cite{BurCLetal20} states that, with probability one, using exact arithmetic,
and in the absence of maximum iteration limits, eventually the convergence
criteria imposed by the parameters $\epsilon_{\ell, \opt}$ and
$\nu_{\ell, \opt}$ must be satisfied.
It is important to note that these stopping criteria, namely
\begin{equation*}
  \begin{aligned}
    \lVert g^{\ell,k_\ell} \rVert_{2} & \leq \nu_{\ell,\opt} & \text{and} &&
    \epsilon_{\ell,k_\ell} & \leq \epsilon_{\ell,\opt}
  \end{aligned}
\end{equation*}
essentially provide an approximate Clarke stationarity certificate.
More precisely, if the parameters $\epsilon_{\ell, \opt}$ and $\nu_{\ell, \opt}$
were set to zero, then all cluster points of the resulting sequence of iterates
must be Clarke stationary for $f^{\ell}$ (see~\cite[Thm.~6.1]{BurCLetal20}),
which amounts to a first-order optimality condition given the Clarke regularity
of~$f^{\ell}$~\cite[p.~753]{BurLO05}.
However, for $\ell < \maxdisc$, no such statement can be made about step $\ell$
in the outer loop of \Cref{alg:amfgs}, because the gradients sampled are not
gradients of $f^{\maxdisc}$.
In contrast, the statement \emph{can} be made about the final step $\ell = L$ in
the outer loop of \Cref{alg:amfgs}.


\section{Numerical experiments}%
\label{sec:examples}

In this section, we present results of applying the new multi-fidelity gradient
sampling algorithms to two applications.
We start by introducing two special cases of the general system~\cref{eqn:sys}
that we will use.
We then describe the experimental setup, and subsequently present the
computational results.


\subsection{Two open-loop systems}%
\label{subsec:two-open-loop}

We test the new methods for the design of $\Hinf$-controllers on two special
instances of open-loop systems~\cref{eqn:sys} that are motivated by
applications discussed subsequently.
First, we consider systems of the form
\begin{equation} \label{eqn:normsys}
  \begin{aligned}
    E \dot{x}(t) & = A x(t) + B w_{1}(t) + B u(t), \\
    z_{1}(t) & = C x(t), \\
    z_{2}(t) & = u(t), \\
    y(t) & = C x(t) + w_{2}(t).
  \end{aligned}
\end{equation}
In~\cref{eqn:normsys}, the disturbances are separated into two
independent parts $w_{1}(t)$ and $w_{2}(t)$, where $w_{1}(t)$ has the same
influence on the system dynamics as the controls and $w_{2}(t)$ disturbs the
measurements taken for the controller.
Also, the performance of the system consists of the non-disturbed
measurements taken for the controller and the control signal itself.
Note that an open-loop system of the form~\cref{eqn:normsys} is known in
the literature as \emph{normalized linear-quadratic Gaussian (LQG)
formulation}; see, e.g.,~\cite{BenHW22, MusG91}.
We may write~\cref{eqn:normsys} in the form~\cref{eqn:sys} by defining
\begin{equation*}
  \begin{aligned}
    B_{1} & = \begin{bmatrix} B & 0 \end{bmatrix}, &
    B_{2} & = B, &
    C_{1} & = \begin{bmatrix} C \\ 0 \end{bmatrix}, &
    C_{2} & = C, \\
    D_{11} & = 0, &
    D_{12} & = \begin{bmatrix} 0 \\ I_{m_{2}} \end{bmatrix}, &
    D_{21} & = \begin{bmatrix} 0 & I_{p_{2}} \end{bmatrix}, &
    D_{22} & = 0.
  \end{aligned}
\end{equation*}

As a second instance of~\cref{eqn:sys}, we consider
\begin{equation} \label{eqn:nnormsys}
  \begin{aligned}
    E \dot{x}(t) & = _{\phantom{1}}A x(t) + _{\phantom{1}}B_{1} w(t) +
      _{\phantom{1}}B_{2} u(t), \\
    z(t) & = C_2 x(t) \phantom{{}+{} _{1}B_{1} w(t)}
      + D_{12} u(t), \\
    y(t) & = C_{2} x(t) + D_{21} w(t).
  \end{aligned}
\end{equation}
Due to the nature of the benchmark problems that we use, the performance and
control measurements are based on the same state observations, i.e.,
we have $C_{1} = C_{2}$ in~\cref{eqn:sys}.
The feed-through term $D_{12}$ is taken as the first
columns ($m_{2} \leq p_{2}$) or rows ($m_{2} > p_{2}$) of the
$\max(m_{2},p_{2})$-dimensional identity matrix, and the feed-through term
$D_{21}$ as the first columns ($m_{1} \leq p_{2}$) or rows ($m_{1} > p_{2}$)
of the $\max(m_{1},p_{2})$-dimensional identity matrix.

For the controller design in both cases, we consider only the problem
formulation of the controller~\cref{eqn:K} without a feed-through term, i.e.,
$\DK = 0$, which is in line with known analytically derived formulas for the
construction of (suboptimal) $\Hinf$-controllers
for~\cref{eqn:normsys,eqn:nnormsys}; see, e.g.,~\cite{BenHW22,DoyGKetal89}.


\subsection{Experimental setup}%
\label{subsec:setup}

\begin{table}[t]
  \caption{Properties of models used in numerical experiments.}
  \label{tab:examples}
  \centering
  {\renewcommand{\arraystretch}{1.25}%
  \settowidth{\maxsteps}{$6\,650$}
  \setlength{\tabcolsep}{.5em}
  \begin{tabular}{rlll}
    \hline
    & & rail example & cylinder example \\
    \hline\noalign{\medskip}
    Discretization levels
      & $\ell = 1$ & $n = 109$ & $n = 6\,618$ \\
    and state dimensions
      & $\ell = 2$ & $n = 371$ & $n = 10\,645$ \\
    & $\ell = 3$ & $n = 1\,357$ & $n = 22\,060$ \\
    & $\ell = 4$ & $n = 5\,177$ & --- \\
    & $\ell = 5$ & $n = 20\,209$ & --- \\
    \noalign{\medskip}\hline\noalign{\medskip}
    Inputs & system~\cref{eqn:normsys} &
      $m_{1} = 13$, $m_{2} = 7$ \hspace{1em}& $m_{1} = 14$, $m_{2} = 6$ \\
    & system~\cref{eqn:nnormsys} &
      $m_{1} = 3$, $m_{2} = 4$ & $m_{1} = 3$, $m_{2} = 3$ \\
    \noalign{\medskip}\hline\noalign{\medskip}
    Outputs & system~\cref{eqn:normsys} &
      $p_{1} = 13$, $p_{2} = 6$ & $p_{1} = 14$, $p_{2} = 8$ \\
      & system~\cref{eqn:nnormsys} &
        $p_{1} = 6$, $p_{2} = 6$ & $p_{1} = 8$, $p_{2} = 8$ \\
    \noalign{\medskip}\hline\noalign{\smallskip}
  \end{tabular}}
\end{table}

We performed our experiments using two publicly available data sets of
spatial discretizations of PDEs~\cite{supWer22b}:
heat flow on a steel bar profile (rail example) and 
laminar fluid flow behind a cylinder obstacle (cylinder example).
The dimensions of the discretizations and the corresponding open-loop systems
are given in \Cref{tab:examples}.
For the cylinder example, the data set provides three different discretizations.
For the rail example, the data set provides nine different discretizations, of
which we chose to use the first five, which allowed us to obtain a sufficiently
accurate approximation while keeping computational costs managable.
We set $\nK$, the order of the controller, to $2$ in all the experiments.

In our experiments, we set the parameters of the multi-fidelity gradient
sampling algorithms as shown in \Cref{tab:parameters}.
While the reduction factors and the line search parameters were set to
default values that do not depend on the discretization level, we chose the
initial sampling radii and stationarity targets to decrease with the increasing
model level.
The rationale for these choices is that the multi-fidelity gradient sampling
algorithms are designed with the idea that final iterates of the optimization on
one level should provide good starting points for the next level, and that as
the level increases it makes sense to set more demanding termination criteria.
Note that we set iteration limits on each level of the multi-fidelity
algorithms.
These values are varied with the problem and are listed in the
column headed ``Max.\ Iters.''\ in the tables that appear below.
In the tables, the point $x^{k_{\ell}}$ denotes the final iterate at level
$\ell$. 
In the case of the rail example, we steadily decrease the maximum number of
allowed iterations per level as the computed iterates approach a minimizer of
the highest fidelity objective.
In the case of the cylinder example, we observed some stagnation in the lowest
fidelity objective for high maximum iteration numbers, perhaps resulting from
a mismatch in the approximation to the highest fidelity objective.
Therefore, we chose here a smaller maximum iteration number than for the
second level.
The number of sampled gradients for all methods and in all problem instances is
set to $q = \nopt +2$, where we recall that $\nopt$, the number of
optimization variables, is given by~\cref{eqn:designvar}. 
The resulting numbers are listed in \Cref{tab:samples}.
All methods are initialized with a randomly generated controller based on 
the same random seed, which is then stabilized by a gradient sampling method
applied to the constraint function~\cref{eqn:hdef}.

\begin{table}[t]
  \caption{Algorithm parameters used in numerical experiments.}
  \label{tab:parameters}
  \centering
  {\renewcommand{\arraystretch}{1.25}%
  \settowidth{\maxsteps}{$6\,650$}
  \setlength{\tabcolsep}{.5em}
  \begin{tabular}{rlll}
    \hline
    & \hfgs{} & \rmfgs{} & \amfgs{} \\
    \hline\noalign{\medskip}
    Init.~sampling radii, &
      $\epsilon_{0} = 0.1$,&
      $\epsilon_{1, 0} = \nu_{1, 0} = 0.1$ &
      $\epsilon_{1, 0} = \nu_{1, 0} = 0.1$ \\
    stationarity targets &  $\nu_{0} = 0.1$ &
      $\epsilon_{2, 0} = \nu_{2, 0} = 0.01$ &
      $\epsilon_{2, 0} = \nu_{2, 0} = 0.01$ \\
    & &
      $\epsilon_{3, 0} = \nu_{3, 0} = 0.001$ &
      $\epsilon_{3, 0} = \nu_{3, 0} = 0.001$ \\
    & &
      $\epsilon_{4, 0} = \nu_{4, 0} = 10^{-4}$ &
      $\epsilon_{4, 0} = \nu_{4, 0} = 10^{-4}$ \\
    & &
      $\epsilon_{5, 0} = \nu_{5, 0} = 10^{-4}$ &
      $\epsilon_{5, 0} = \nu_{5, 0} = 10^{-4}$ \\
    \noalign{\medskip}\hline\noalign{\medskip}
    Termination tol. &
      $\epsilon_{\opt} = 10^{-4}$,  &
      $\epsilon_{1, \opt} = \nu_{1, \opt} = 10^{-4}$ &
      $\epsilon_{1, \opt} = \nu_{1, \opt} = 0.01$ \\
    & $\nu_{\opt} =10^{-4}$& 
      $\epsilon_{2, \opt} = \nu_{2, \opt} = 10^{-4}$ &
      $\epsilon_{2, \opt} = \nu_{2, \opt} = 0.001$ \\
    & & 
      $\epsilon_{3, \opt} = \nu_{3, \opt} = 10^{-4}$ &
      $\epsilon_{3, \opt} = \nu_{3, \opt} = 10^{-4}$ \\
    & & 
      $\epsilon_{4, \opt} = \nu_{4, \opt} = 10^{-4}$ &
      $\epsilon_{4, \opt} = \nu_{4, \opt} = 10^{-4}$ \\
    & & 
      $\epsilon_{5, \opt} = \nu_{5, \opt} = 10^{-4}$ &
      $\epsilon_{5, \opt} = \nu_{5, \opt} = 10^{-4}$ \\
    \noalign{\medskip}\hline\noalign{\medskip}
    Reduction factors &
      $\theta_{\epsilon} = 0.1$, &
      $\theta_{\ell, \epsilon} = \theta_{\ell, \nu} = 0.1$ &
      $\theta_{\ell, \epsilon} = \theta_{\ell, \nu} = 0.1$ \\
    & $\theta_{\nu} = 0.1$ &
      for $\ell = 1, \ldots, \maxdisc$ &
      for $\ell = 1, \ldots, \maxdisc$ \\
    \noalign{\medskip}\hline\noalign{\medskip}
    Line search &
      $\beta = 10^{-4}$, &
      $\beta_{\ell} = 10^{-4}$ &
      $\beta_{\ell} = 10^{-4}$ \\
    &
      $\gamma = 0.5$ &
      $\gamma_{\ell} = 0.5$ &
      $\gamma_{\ell} = 0.5$ \\
    & &
      for $\ell = 1, \ldots, \maxdisc$ &
      for $\ell = 1, \ldots, \maxdisc$ \\
    \noalign{\medskip}\hline\noalign{\smallskip}
  \end{tabular}}
\end{table}

\begin{table}[t]
  \caption{Number of sampled gradients per problem instance.}
  \label{tab:samples}
  \centering
  {\renewcommand{\arraystretch}{1.25}%
  \settowidth{\maxsteps}{$6\,650$}
  \setlength{\tabcolsep}{.5em}
  \begin{tabular}{lcccc}
    \hline
    & \multicolumn{2}{c}{rail example} & \multicolumn{2}{c}{cylinder example} \\
    & system~\cref{eqn:normsys} & system~\cref{eqn:nnormsys} &
      system~\cref{eqn:normsys} & system~\cref{eqn:nnormsys} \\
    \hline\noalign{\medskip}
    \# sampled gradients $q$ &
      $32$ &
      $26$ &
      $34$ &
      $24$ \\
    \noalign{\medskip}\hline\noalign{\smallskip}
  \end{tabular}}
\end{table}

We compare \rmfgs{} and \amfgs{} to the
single-fidelity gradient sampling method from~\cite[Alg.~GS]{BurCLetal20}
applied directly to the high-fidelity objective function $f^{\maxdisc}$,
denoted subsequently as \hfgs{}.
We compare the results for the different methods by comparing
the evolution of the high-fidelity objective $f^{\maxdisc}$ on the iterate
sequence $\{x^{k}\}$.
In the case of \rmfgs{}, which does not access $f^{\maxdisc}$ until its final
outer loop, we computed $f^{\maxdisc}(x^{k})$ \emph{a posteriori}.

For each problem instance that we solve, since we do not know the minimal value
of $f^{\maxdisc}$, it is convenient to define
\begin{equation*}
  f_{\min} := \min\left(
    f^{\maxdisc}(x_{\mathrm{\hfgs}}),
    f^{\maxdisc}(x_{\mathrm{\rmfgs}}),
    f^{\maxdisc}(x_{\mathrm{\amfgs}})
    \right),
\end{equation*}
where the three quantities on the right-hand side are respectively the minimal
values of $f^{\maxdisc}$ found by the three different methods.
Then, in the figures below, for each problem instance we show two different
plots of the evolution of $f^{\maxdisc}(x^{k})$.
In the plots on the left, the vertical axis shows the values of $f^{\maxdisc}$
computed by each of the three methods, with different symbols indicating the
discretization level, i.e., the index of the outer loop in the case of \rmfgs{}
and \amfgs{}.
For \hfgs{}, only the highest fidelity discretization symbol is used.
In the plots on the right, the vertical axis shows the relative error
\begin{equation*}
  \frac{f^{\maxdisc}(x^{k}) - f_{\min}}{f_{\min}},
\end{equation*}
using $f_{\min}$ as our best estimate of the true minimal value.
In both cases, the horizontal axis shows the running time in hours.

The experiments were run on compute nodes of the \texttt{Greene}
high-performance computing cluster of the New York University using 16
processing cores of the Intel Xeon Platinum 8268 24C 205W CPU at 2.90\,GHz
and 16\,GB main memory.
We used MATLAB 9.9.0.1467703 (R2020b) running on Red Hat Enterprise Linux
release 8.4 (Ootpa).
For the single-fidelity gradient sampling method, we used the implementation
in HANSO, Hybrid Algorithm for Non-Smooth Optimization,
version 3.0~\cite{Ove21}.
The new multi-fidelity codes are also based on this.
All the examples discussed below, except the first
two levels of the rail example, use MATLAB's sparse data structure.
For the computation of the $\Hinf$-norm we employ the \texttt{normTfMaxPeak} and
\texttt{normTfPeak} routines from ROSTAPACK (RObust STAbility PACKage),
version 3.0~\cite{Mit22}; see
also~\cite{BenM18} for the implemented algorithms.
As \texttt{normTfPeak} does not do a stability check,  we implemented this using
MATLAB's \texttt{eigs} function.
The source code, data and results of the numerical experiments are
open source/open access and available at~\cite{supWer22b}.


\subsection{Optimal cooling of a steel rail profile}

We consider the heat flow on a two-dimensional cross section of a steel bar
for optimal cooling; see~\cite{Saa09} for further details and~\cite{SaaKB21}
for the data set.
The underlying heat equation is discretized on multiple grid levels using finite
elements.
The resulting dimensions of the two open-loop
systems~\cref{eqn:normsys,eqn:nnormsys} can be found in the rail
example column of \Cref{tab:examples}.

\begin{figure}[t]
  \centering
  \begin{subfigure}[b]{.49\linewidth}
    \centering
  \tikzexternalenable%
  \tikzsetnextfilename{rail_norm_objtime}%
  \begin{tikzpicture}[font = \small]
  \begin{axis}[%
    name   = states,
    width  = .7\textwidth,
    height = .18\textheight,
    scale only axis,
    xmin = 0,
    xmax = 55,
    ymin = 0.435,
    ymax = 0.5,
    xminorticks = false,
    yminorticks = false,
    xlabel = {Time (h)},
    ylabel = {Objective function $f^{\maxdisc}(x^{k})$},
    ylabel style   = {yshift = -.3em},
    scaled x ticks = false,
    x tick label style = {/pgf/number format/1000 sep={\,}},
    y tick label style = {/pgf/number format/1000 sep={\,}},
    cycle list name    = railplotlist
  ]
    \pgfplotstableread{graphics/data/rail_norm_objtime_gs.dat}\tableDIRECT
    \addplot+[mark repeat = 17] table[x index = 0, y index = 1] {\tableDIRECT};
    
    \pgfplotstableread{graphics/data/rail_norm_objtime_rmfgs1.dat}\tableMULTI
    \addplot+[mark repeat = 34] table[x index = 0, y index = 1] {\tableMULTI};
    
    \pgfplotstableread{graphics/data/rail_norm_objtime_rmfgs2.dat}\tableMULTI
    \addplot+[mark repeat = 34] table[x index = 0, y index = 1] {\tableMULTI};
    
    \pgfplotstableread{graphics/data/rail_norm_objtime_rmfgs3.dat}\tableMULTI
    \addplot+[mark repeat = 34] table[x index = 0, y index = 1] {\tableMULTI};
    
    \pgfplotstableread{graphics/data/rail_norm_objtime_rmfgs4.dat}\tableMULTI
    \addplot+[mark repeat = 34] table[x index = 0, y index = 1] {\tableMULTI};
    
    \pgfplotstableread{graphics/data/rail_norm_objtime_rmfgs5.dat}\tableMULTI
    \addplot+[mark repeat = 17] table[x index = 0, y index = 1] {\tableMULTI};
    
    \pgfplotstableread{graphics/data/rail_norm_objtime_amfgs1.dat}\tableMULTI
    \addplot+[mark repeat = 22] table[x index = 0, y index = 1] {\tableMULTI};
    
    \pgfplotstableread{graphics/data/rail_norm_objtime_amfgs2.dat}\tableMULTI
    \addplot+[mark repeat = 334] table[x index = 0, y index = 1] {\tableMULTI};
    
    \pgfplotstableread{graphics/data/rail_norm_objtime_amfgs3.dat}\tableMULTI
    \addplot+[mark repeat = 167] table[x index = 0, y index = 1] {\tableMULTI};
    
    \pgfplotstableread{graphics/data/rail_norm_objtime_amfgs4.dat}\tableMULTI
    \addplot+[mark repeat = 34] table[x index = 0, y index = 1] {\tableMULTI};
    
    \pgfplotstableread{graphics/data/rail_norm_objtime_amfgs5.dat}\tableMULTI
    \addplot+[mark repeat = 17] table[x index = 0, y index = 1] {\tableMULTI};
  \end{axis}
\end{tikzpicture}%
  \tikzexternaldisable%

    \caption{objective function values}
    \label{fig:rail_norm_objtime}
  \end{subfigure}%
  \hfill%
  \begin{subfigure}[b]{.49\linewidth}
    \centering
  \tikzexternalenable%
  \tikzsetnextfilename{rail_norm_acctime}%
  \begin{tikzpicture}[font = \small]
  \begin{semilogyaxis}[%
    name   = states,
    width  = .7\textwidth,
    height = .18\textheight,
    scale only axis,
    xmin = 0,
    xmax = 55,
    ymin = 1e-6,
    ymax = 2e-01,
    xminorticks = false,
    yminorticks = false,
    xlabel = {Time (h)},
    ylabel = {$(f^{\maxdisc}(x^{k}) - f_{\min}) / f_{\min}$},
    ylabel style   = {yshift = -.3em},
    scaled x ticks = false,
    x tick label style = {/pgf/number format/1000 sep={\,}},
    y tick label style = {/pgf/number format/1000 sep={\,}},
    cycle list name    = railplotlist
  ]
    \pgfplotstableread{graphics/data/rail_norm_acctime_gs.dat}\tableDIRECT
    \addplot+[mark repeat = 17] table[x index = 0, y index = 1] {\tableDIRECT};
    
    \pgfplotstableread{graphics/data/rail_norm_acctime_rmfgs1.dat}\tableMULTI
    \addplot+[mark repeat = 34] table[x index = 0, y index = 1] {\tableMULTI};
    
    \pgfplotstableread{graphics/data/rail_norm_acctime_rmfgs2.dat}\tableMULTI
    \addplot+[mark repeat = 34] table[x index = 0, y index = 1] {\tableMULTI};
    
    \pgfplotstableread{graphics/data/rail_norm_acctime_rmfgs3.dat}\tableMULTI
    \addplot+[mark repeat = 34] table[x index = 0, y index = 1] {\tableMULTI};
    
    \pgfplotstableread{graphics/data/rail_norm_acctime_rmfgs4.dat}\tableMULTI
    \addplot+[mark repeat = 34] table[x index = 0, y index = 1] {\tableMULTI};
    
    \pgfplotstableread{graphics/data/rail_norm_acctime_rmfgs5.dat}\tableMULTI
    \addplot+[mark repeat = 17] table[x index = 0, y index = 1] {\tableMULTI};
    
    \pgfplotstableread{graphics/data/rail_norm_acctime_amfgs1.dat}\tableMULTI
    \addplot+[mark repeat = 22] table[x index = 0, y index = 1] {\tableMULTI};
    
    \pgfplotstableread{graphics/data/rail_norm_acctime_amfgs2.dat}\tableMULTI
    \addplot+[mark repeat = 334] table[x index = 0, y index = 1] {\tableMULTI};
    
    \pgfplotstableread{graphics/data/rail_norm_acctime_amfgs3.dat}\tableMULTI
    \addplot+[mark repeat = 167] table[x index = 0, y index = 1] {\tableMULTI};
    
    \pgfplotstableread{graphics/data/rail_norm_acctime_amfgs4.dat}\tableMULTI
    \addplot+[mark repeat = 34] table[x index = 0, y index = 1] {\tableMULTI};
    
    \pgfplotstableread{graphics/data/rail_norm_acctime_amfgs5.dat}\tableMULTI
    \addplot+[mark repeat = 17] table[x index = 0, y index = 1] {\tableMULTI};
  \end{semilogyaxis}
\end{tikzpicture}%
  \tikzexternaldisable%

    \caption{distance to lowest objective value}
    \label{fig:rail_norm_acctime}
  \end{subfigure}
  
  \vspace{.5\baselineskip}
  \tikzexternalenable%
  \tikzsetnextfilename{rail_legend}%
  \begin{tikzpicture}[font = \small]
  \begin{axis}[%
    hide axis,
    width  = .7\textwidth,
    height = .1\textheight,
    scale only axis,
    xmin = 0,
    xmax = 1,
    ymin = 0,
    ymax = 1,
    legend columns = 5,
    legend cell align = {left},
    legend style = {
      at     = {(0,0)},
      anchor = center,
      /tikz/every even column/.append style = {column sep = 0.5cm}}
  ]
    
    \addlegendimage{gs}
    \addlegendentry{\hfgs{}}
    
    \addlegendimage{rmfgs}
    \addlegendentry{\rmfgs{}}
    
    \addlegendimage{amfgs}
    \addlegendentry{\amfgs{}}
    
    \pgfplotsinvokeforeach{2, 3}{
      \addlegendimage{amfgs,white}
      \addlegendentry{}
    }
    
    \pgfplotsset{cycle list name = railplotlist, cycle list shift = 1}
    \pgfplotsinvokeforeach{1, 2, ..., 5}{
      \addplot+[only marks, black] coordinates {(0,0)};
    }
    \addlegendentry{level $1$}
    \addlegendentry{level $2$}
    \addlegendentry{level $3$}
    \addlegendentry{level $4$}
    \addlegendentry{level $5$}
  \end{axis}
\end{tikzpicture}%
  \tikzexternaldisable%

  \caption{Rail example with formulation~\cref{eqn:normsys}:
    To reach the final objective function value found by \hfgs{}, \rmfgs{}
    achieves a speedup of $452$ and \amfgs{} achieves a speedup of $30$ in
    comparison.
    Additionally, \rmfgs{} and \amfgs{} ultimately obtain 
    lower objective function values than those found by using only the
    high-fidelity model in \hfgs{}.}
  \label{fig:rail_norm}
\end{figure}
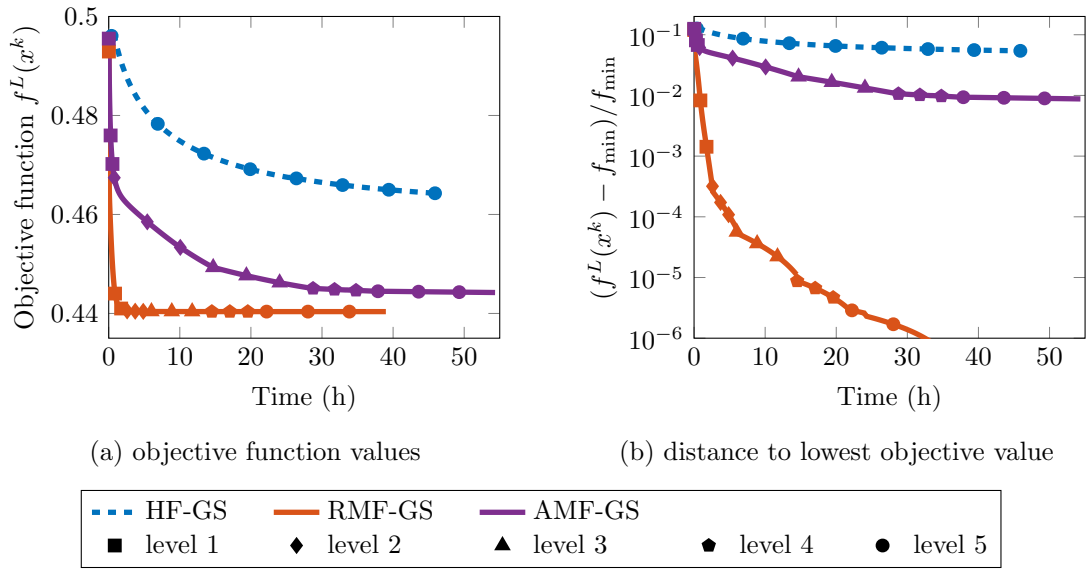

\begin{table}[t]
  \caption{Rail example with formulation~\cref{eqn:normsys}:
    The table reports the wall-clock time of the computations, the
    number of iterations taken versus the maximum allowed number and
    the objective function values corresponding to the low-fidelity models
    (in case of \rmfgs{}) and high-fidelity models.}
  \label{tab:rail_norm}
  \centering
  {\renewcommand{\arraystretch}{1.25}%
  \settowidth{\maxsteps}{$6\,650$}
  \setlength{\tabcolsep}{.5em}
  \begin{tabular}{lrrrrr}
    \hline
    & &
      Time (h) &
      Iters./Max.\ Iters. &
      $f^{\ell}(x^{k_{\ell}})$ &
      $f^{\maxdisc}(x^{k_{\ell}})$ \\
    \hline\noalign{\medskip}
    \hfgs{} &
      & $45.895$
      & $120$~/~\makebox[\maxsteps][r]{$120$}
      & ---
      & $0.464284$ \\
    \noalign{\medskip}\hline\noalign{\medskip}
    \rmfgs{} & level $1$
      & $2.5594$
      & $5\,000$~/~\makebox[\maxsteps][r]{$5\,000$}
      & $0.440143$
      & $0.440511$ \\
    & level $2$
      & $3.3656$
      & $1\,000$~/~\makebox[\maxsteps][r]{$1\,000$}
      & $0.440312$
      & $0.440399$ \\
    & level $3$
      & $8.5062$
      & $500$~/~\makebox[\maxsteps][r]{$500$}
      & $0.440365$
      & $0.440375$ \\
    & level $4$
      & $7.4379$
      & $100$~/~\makebox[\maxsteps][r]{$100$}
      & $0.440372$
      & $0.440372$ \\
    & level $5$
      & $17.113$
      & $50$~/~\makebox[\maxsteps][r]{$50$}
      & ---
      & $0.440370$ \\
    \noalign{\smallskip}\cline{2-6}\noalign{\smallskip}
    & 
      & $38.982$
      & $6\,650$~/~\makebox[\maxsteps][r]{$6\,650$}
      & ---
      & $0.440370$ \\
    \noalign{\smallskip}\hline\noalign{\medskip}
    \amfgs{} & level $1$
      & $0.7683$
      & $66$~/~\makebox[\maxsteps][r]{$5\,000$}
      & ---
      & $0.467422$ \\
    & level $2$
      & $13.901$
      & $1\,000$~/~\makebox[\maxsteps][r]{$1\,000$}
      & ---
      & $0.449315$ \\
    & level $3$
      & $13.983$
      & $500$~/~\makebox[\maxsteps][r]{$500$}
      & ---
      & $0.445053$ \\
    & level $4$
      & $8.8870$
      & $100$~/~\makebox[\maxsteps][r]{$100$}
      & ---
      & $0.444489$ \\
    & level $5$
      & $16.805$
      & $50$~/~\makebox[\maxsteps][r]{$50$}
      & ---
      & $0.444229$ \\
    \noalign{\smallskip}\cline{2-6}\noalign{\smallskip}
    & 
      & $54.363$
      & $1\,716$~/~\makebox[\maxsteps][r]{$6\,650$}
      & ---
      & $0.444229$ \\
    \noalign{\smallskip}\hline\noalign{\smallskip}
  \end{tabular}}
\end{table}
 
We first consider the example formulation~\cref{eqn:normsys}.
The results are shown in \Cref{fig:rail_norm,tab:rail_norm}.
Even a quick glance reveals that both new methods are faster
and more accurate than the single-fidelity method \hfgs{}, with \rmfgs{}
faster and more accurate than \amfgs{}.
Indeed, already level $1$ of the \rmfgs{} method obtains in less than
$0.1$\,h about the same value for $f^{\maxdisc}$ as the final value found by
\hfgs{} after $45$\,h.
Furthermore, although the plot on the left side of \Cref{fig:rail_norm} suggests
that \rmfgs{} stagnates, the plot on the right side shows that this
is not the case, with additional digits of accuracy steadily attained as the
hierarchy level of \rmfgs{} is increased.
Overall, \rmfgs{} achieves a speedup of $452$ compared to \hfgs{} to reach
the same high-fidelity objective function value.
\amfgs{} achieves a speedup of $30$ compared to \hfgs{}.
For all methods, the stabilization of the initial guess took only a
single step of gradient sampling for the spectral abscissa constraint function.
Even for \Cref{alg:rmfgs}, no subsequent stabilization steps were required.

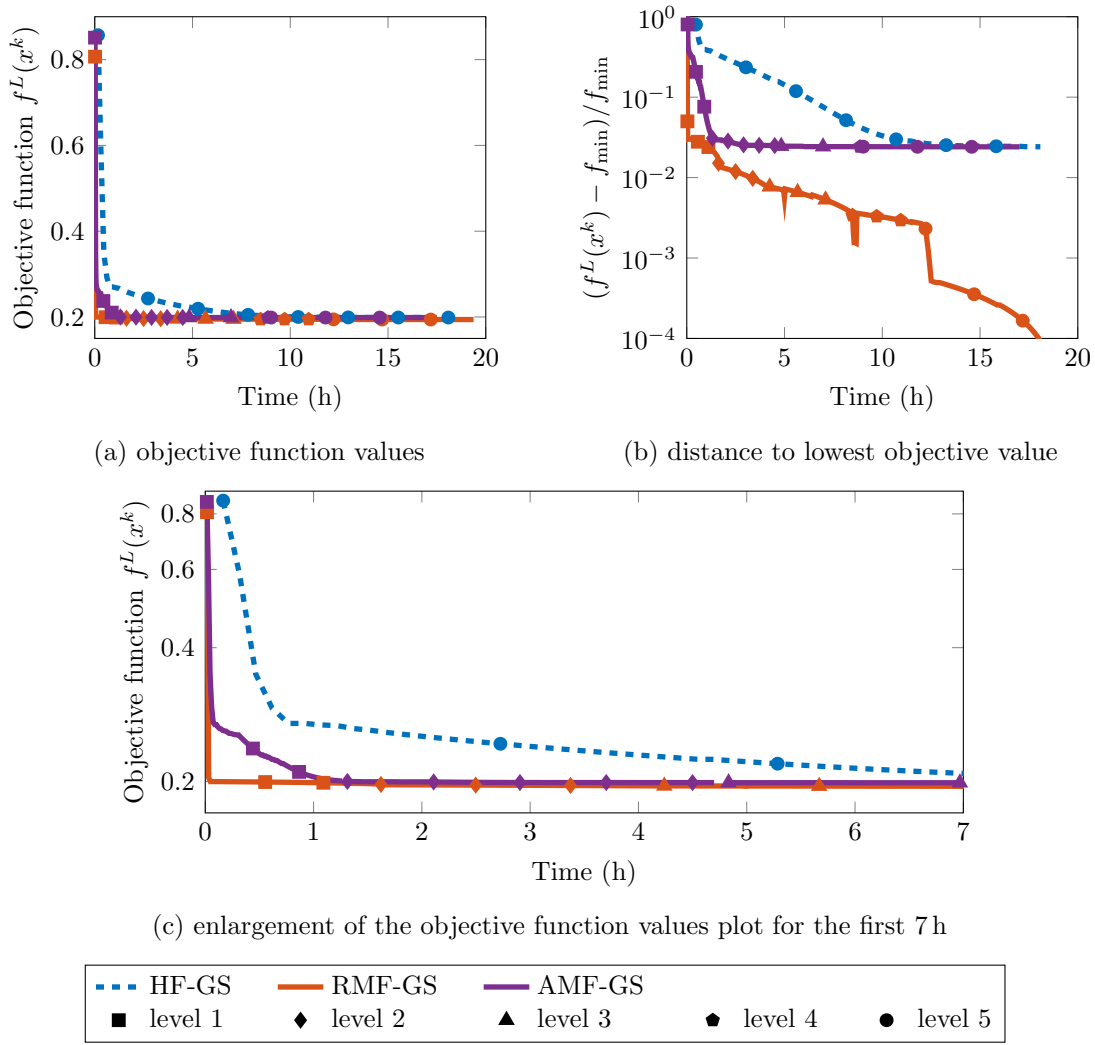
\begin{figure}[t]
  \centering
  \begin{subfigure}[b]{.49\linewidth}
    \centering
  \tikzexternalenable%
  \tikzsetnextfilename{rail_objtime}%
  \begin{tikzpicture}[font = \small]
  \begin{axis}[%
    name   = states,
    width  = .7\textwidth,
    height = .18\textheight,
    scale only axis,
    xmin = 0,
    xmax = 20,
    ymin = 0.15,
    ymax = 0.9,
    xminorticks = false,
    yminorticks = false,
    xlabel = {Time (h)},
    ylabel = {Objective function $f^{\maxdisc}(x^{k})$},
    ylabel style   = {yshift = -.3em},
    scaled x ticks = false,
    x tick label style = {/pgf/number format/1000 sep={\,}},
    y tick label style = {/pgf/number format/1000 sep={\,}},
    cycle list name    = railplotlist
  ]
    \pgfplotstableread{graphics/data/rail_objtime_gs.dat}\tableDIRECT
    \addplot+[mark repeat = 17] table[x index = 0, y index = 1] {\tableDIRECT};
    
    \pgfplotstableread{graphics/data/rail_objtime_rmfgs1.dat}\tableMULTI
    \addplot+[mark repeat = 34] table[x index = 0, y index = 1] {\tableMULTI};
    
    \pgfplotstableread{graphics/data/rail_objtime_rmfgs2.dat}\tableMULTI
    \addplot+[mark repeat = 34] table[x index = 0, y index = 1] {\tableMULTI};
    
    \pgfplotstableread{graphics/data/rail_objtime_rmfgs3.dat}\tableMULTI
    \addplot+[mark repeat = 34] table[x index = 0, y index = 1] {\tableMULTI};
    
    \pgfplotstableread{graphics/data/rail_objtime_rmfgs4.dat}\tableMULTI
    \addplot+[mark repeat = 34] table[x index = 0, y index = 1] {\tableMULTI};
    
    \pgfplotstableread{graphics/data/rail_objtime_rmfgs5.dat}\tableMULTI
    \addplot+[mark repeat = 17] table[x index = 0, y index = 1] {\tableMULTI};
    
    \pgfplotstableread{graphics/data/rail_objtime_amfgs1.dat}\tableMULTI
    \addplot+[mark repeat = 29] table[x index = 0, y index = 1] {\tableMULTI};
    
    \pgfplotstableread{graphics/data/rail_objtime_amfgs2.dat}\tableMULTI
    \addplot+[mark repeat = 16] table[x index = 0, y index = 1] {\tableMULTI};
    
    \pgfplotstableread{graphics/data/rail_objtime_amfgs3.dat}\tableMULTI
    \addplot+[mark repeat = 16] table[x index = 0, y index = 1] {\tableMULTI};
    
    \pgfplotstableread{graphics/data/rail_objtime_amfgs4.dat}\tableMULTI
    \addplot+[mark repeat = 1] table[x index = 0, y index = 1] {\tableMULTI};
    
    \pgfplotstableread{graphics/data/rail_objtime_amfgs5.dat}\tableMULTI
    \addplot+[mark repeat = 17] table[x index = 0, y index = 1] {\tableMULTI};
  \end{axis}
\end{tikzpicture}%
  \tikzexternaldisable%

    \caption{objective function values}
    \label{fig:rail_objtime}
  \end{subfigure}%
  \hfill%
  \begin{subfigure}[b]{.49\linewidth}
    \centering
  \tikzexternalenable%
  \tikzsetnextfilename{rail_acctime}%
  \begin{tikzpicture}[font = \small]
  \begin{semilogyaxis}[%
    name   = states,
    width  = .7\textwidth,
    height = .18\textheight,
    scale only axis,
    xmin = 0,
    xmax = 20,
    ymin = 1e-04,
    ymax = 1e+00,
    xminorticks = false,
    yminorticks = false,
    xlabel = {Time (h)},
    ylabel = {$(f^{\maxdisc}(x^{k}) - f_{\min}) / f_{\min}$},
    ylabel style   = {yshift = -.3em},
    scaled x ticks = false,
    x tick label style = {/pgf/number format/1000 sep={\,}},
    y tick label style = {/pgf/number format/1000 sep={\,}},
    cycle list name    = railplotlist
  ]
    \pgfplotstableread{graphics/data/rail_acctime_gs.dat}\tableDIRECT
    \addplot+[mark repeat = 17] table[x index = 0, y index = 1] {\tableDIRECT};
    
    \pgfplotstableread{graphics/data/rail_acctime_rmfgs1.dat}\tableMULTI
    \addplot+[mark repeat = 34] table[x index = 0, y index = 1] {\tableMULTI};
    
    \pgfplotstableread{graphics/data/rail_acctime_rmfgs2.dat}\tableMULTI
    \addplot+[mark repeat = 34] table[x index = 0, y index = 1] {\tableMULTI};
    
    \pgfplotstableread{graphics/data/rail_acctime_rmfgs3.dat}\tableMULTI
    \addplot+[mark repeat = 34] table[x index = 0, y index = 1] {\tableMULTI};
    
    \pgfplotstableread{graphics/data/rail_acctime_rmfgs4.dat}\tableMULTI
    \addplot+[mark repeat = 34] table[x index = 0, y index = 1] {\tableMULTI};
    
    \pgfplotstableread{graphics/data/rail_acctime_rmfgs5.dat}\tableMULTI
    \addplot+[mark repeat = 17] table[x index = 0, y index = 1] {\tableMULTI};
    
    \pgfplotstableread{graphics/data/rail_acctime_amfgs1.dat}\tableMULTI
    \addplot+[mark repeat = 29] table[x index = 0, y index = 1] {\tableMULTI};
    
    \pgfplotstableread{graphics/data/rail_acctime_amfgs2.dat}\tableMULTI
    \addplot+[mark repeat = 16] table[x index = 0, y index = 1] {\tableMULTI};
    
    \pgfplotstableread{graphics/data/rail_acctime_amfgs3.dat}\tableMULTI
    \addplot+[mark repeat = 16] table[x index = 0, y index = 1] {\tableMULTI};
    
    \pgfplotstableread{graphics/data/rail_acctime_amfgs4.dat}\tableMULTI
    \addplot+[mark repeat = 1] table[x index = 0, y index = 1] {\tableMULTI};
    
    \pgfplotstableread{graphics/data/rail_acctime_amfgs5.dat}\tableMULTI
    \addplot+[mark repeat = 17] table[x index = 0, y index = 1] {\tableMULTI};
  \end{semilogyaxis}
\end{tikzpicture}%
  \tikzexternaldisable%

    \caption{distance to lowest objective value}
    \label{fig:rail_acctime}
  \end{subfigure}%
  
  \vspace{.5\baselineskip}
  \begin{subfigure}[b]{.95\linewidth}
    \centering
  \tikzexternalenable%
  \tikzsetnextfilename{rail_objtime_zoom}%
  \begin{tikzpicture}[font = \small]
  \begin{semilogyaxis}[%
    name   = states,
    width  = .7\textwidth,
    height = .18\textheight,
    scale only axis,
    xmin = 0,
    xmax = 7,
    ymin = 0.17,
    ymax = 0.9,
    xminorticks = false,
    yminorticks = false,
    ytick       = {0.8, 0.6, 0.4, 0.2},
    yticklabels = {0.8, 0.6, 0.4, 0.2},
    xlabel = {Time (h)},
    ylabel = {Objective function $f^{\maxdisc}(x^{k})$},
    ylabel style   = {yshift = -.3em},
    scaled x ticks = false,
    x tick label style = {/pgf/number format/1000 sep={\,}},
    y tick label style = {/pgf/number format/1000 sep={\,}},
    cycle list name    = railplotlist
  ]
    \pgfplotstableread{graphics/data/rail_objtime_gs.dat}\tableDIRECT
    \addplot+[mark repeat = 17] table[x index = 0, y index = 1] {\tableDIRECT};
    
    \pgfplotstableread{graphics/data/rail_objtime_rmfgs1.dat}\tableMULTI
    \addplot+[mark repeat = 34] table[x index = 0, y index = 1] {\tableMULTI};
    
    \pgfplotstableread{graphics/data/rail_objtime_rmfgs2.dat}\tableMULTI
    \addplot+[mark repeat = 34] table[x index = 0, y index = 1] {\tableMULTI};
    
    \pgfplotstableread{graphics/data/rail_objtime_rmfgs3.dat}\tableMULTI
    \addplot+[mark repeat = 34] table[x index = 0, y index = 1] {\tableMULTI};
    
    \pgfplotstableread{graphics/data/rail_objtime_rmfgs4.dat}\tableMULTI
    \addplot+[mark repeat = 34] table[x index = 0, y index = 1] {\tableMULTI};
    
    \pgfplotstableread{graphics/data/rail_objtime_rmfgs5.dat}\tableMULTI
    \addplot+[mark repeat = 17] table[x index = 0, y index = 1] {\tableMULTI};
    
    \pgfplotstableread{graphics/data/rail_objtime_amfgs1.dat}\tableMULTI
    \addplot+[mark repeat = 29] table[x index = 0, y index = 1] {\tableMULTI};
    
    \pgfplotstableread{graphics/data/rail_objtime_amfgs2.dat}\tableMULTI
    \addplot+[mark repeat = 16] table[x index = 0, y index = 1] {\tableMULTI};
    
    \pgfplotstableread{graphics/data/rail_objtime_amfgs3.dat}\tableMULTI
    \addplot+[mark repeat = 16] table[x index = 0, y index = 1] {\tableMULTI};
    
    \pgfplotstableread{graphics/data/rail_objtime_amfgs4.dat}\tableMULTI
    \addplot+[mark repeat = 1] table[x index = 0, y index = 1] {\tableMULTI};
    
    \pgfplotstableread{graphics/data/rail_objtime_amfgs5.dat}\tableMULTI
    \addplot+[mark repeat = 17] table[x index = 0, y index = 1] {\tableMULTI};
  \end{semilogyaxis}
\end{tikzpicture}%
  \tikzexternaldisable%

    \caption{enlargement of the objective function values plot for the first
      $7$\,h}
    \label{fig:rail_objtime_zoom}
  \end{subfigure}%

  \vspace{.5\baselineskip}
  \tikzexternalenable%
  \tikzsetnextfilename{rail_legend}%
  \begin{tikzpicture}[font = \small]
  \begin{axis}[%
    hide axis,
    width  = .7\textwidth,
    height = .1\textheight,
    scale only axis,
    xmin = 0,
    xmax = 1,
    ymin = 0,
    ymax = 1,
    legend columns = 5,
    legend cell align = {left},
    legend style = {
      at     = {(0,0)},
      anchor = center,
      /tikz/every even column/.append style = {column sep = 0.5cm}}
  ]
    
    \addlegendimage{gs}
    \addlegendentry{\hfgs{}}
    
    \addlegendimage{rmfgs}
    \addlegendentry{\rmfgs{}}
    
    \addlegendimage{amfgs}
    \addlegendentry{\amfgs{}}
    
    \pgfplotsinvokeforeach{2, 3}{
      \addlegendimage{amfgs,white}
      \addlegendentry{}
    }
    
    \pgfplotsset{cycle list name = railplotlist, cycle list shift = 1}
    \pgfplotsinvokeforeach{1, 2, ..., 5}{
      \addplot+[only marks, black] coordinates {(0,0)};
    }
    \addlegendentry{level $1$}
    \addlegendentry{level $2$}
    \addlegendentry{level $3$}
    \addlegendentry{level $4$}
    \addlegendentry{level $5$}
  \end{axis}
\end{tikzpicture}%
  \tikzexternaldisable%

  \caption{Rail example with formulation~\cref{eqn:nnormsys}:
    To reach the final objective function value found by \hfgs{}, \rmfgs{}
    achieves a speedup of $17$ and \amfgs{} a speedup of $2$ in comparison.}
  \label{fig:rail}
\end{figure}

\begin{table}[t]
  \caption{Rail example with formulation~\cref{eqn:nnormsys}:
    The table reports the wall-clock time of the computations, the
    number of iterations taken versus the maximum allowed number and
    the objective function values corresponding to the low-fidelity models (in
    case of \rmfgs{}) and high-fidelity models.}
  \label{tab:rail}
  \centering
  {\renewcommand{\arraystretch}{1.25}%
  \settowidth{\maxsteps}{$6\,650$}
  \setlength{\tabcolsep}{.5em}
  \begin{tabular}{lrrrrr}
    \hline
    & &
      Time (h) &
      Iters./Max.\ Iters. &
      $f^{\ell}(x^{k_{\ell}})$ &
      $f^{\maxdisc}(x^{k_{\ell}})$ \\
    \hline\noalign{\medskip}
    \hfgs{} &
      & $18.085$
      & $120$~/~\makebox[\maxsteps][r]{$120$}
      & ---
      & $0.198720$ \\
    \noalign{\medskip}\hline\noalign{\medskip}
    \rmfgs{} & level $1$
      & $1.5963$
      & $5\,000$~/~\makebox[\maxsteps][r]{$5\,000$}
      & $0.197222$
      & $0.197473$ \\
    & level $2$
      & $2.6002$
      & $1\,000$~/~\makebox[\maxsteps][r]{$1\,000$}
      & $0.195428$
      & $0.195475$ \\
    & level $3$
      & $4.2480$
      & $500$~/~\makebox[\maxsteps][r]{$500$}
      & $0.194404$
      & $0.194740$ \\
    & level $4$
      & $3.6196$
      & $100$~/~\makebox[\maxsteps][r]{$100$}
      & $0.194148$
      & $0.194543$ \\
    & level $5$
      & $7.3114$
      & $50$~/~\makebox[\maxsteps][r]{$50$}
      & ---
      & $0.194028$ \\
    \noalign{\smallskip}\cline{2-6}\noalign{\smallskip}
    & 
      & $19.375$
      & $6\,650$~/~\makebox[\maxsteps][r]{$6\,650$}
      & ---
      & $0.194028$ \\
    \noalign{\smallskip}\hline\noalign{\medskip}
    \amfgs{} & level $1$
      & $1.2626$
      & $86$~/~\makebox[\maxsteps][r]{$5\,000$}
      & ---
      & $0.200531$ \\
    & level $2$
      & $3.4359$
      & $69$~/~\makebox[\maxsteps][r]{$1\,000$}
      & ---
      & $0.198870$ \\
    & level $3$
      & $3.8725$
      & $29$~/~\makebox[\maxsteps][r]{$500$}
      & ---
      & $0.198732$ \\
    & level $4$
      & $0.2885$
      & $1$~/~\makebox[\maxsteps][r]{$100$}
      & ---
      & $0.198732$ \\
    & level $5$
      & $8.1653$
      & $50$~/~\makebox[\maxsteps][r]{$50$}
      & ---
      & $0.198707$ \\
    \noalign{\smallskip}\cline{2-6}\noalign{\smallskip}
    & 
      & $17.043$
      & $235$~/~\makebox[\maxsteps][r]{$6\,650$}
      & ---
      & $0.198707$ \\
    \noalign{\smallskip}\hline\noalign{\smallskip}
  \end{tabular}}
\end{table}

The second experiment that we consider for this application is for
formulation~\cref{eqn:nnormsys}.
The disturbances are set to be the lower boundary temperatures and the
controls are restricted to the boundary temperatures of the upper segments;
see also~\cite[Sec.~3.2]{BenHW22a} where the same setup is used.
The results are shown in \Cref{fig:rail,tab:rail}.
In this case, although the results in absolute terms are not as
much in favor of the new methods as they were for the previous
example, in relative terms, \rmfgs{} is much better than either of the other
methods, and \amfgs{} gives much better results than the single-fidelity method
until after $10$\,h of computation.
\rmfgs{} and \amfgs{} reach the same level of the final objective function
value of \hfgs{} in about $1.5$\,h and both provide at the end of the
iterations a smaller objective function value than \hfgs{}.
All methods needed only a single gradient sampling step to stabilize the
closed-loop system at initialization.


\subsection{Robust stabilization of laminar flows in a cylinder wake}

We now consider the stabilization of laminar flow in a two-dimensional wake  
resulting from a circular obstacle.
The flow is modeled as the linearization of the Navier-Stokes equations at
Reynolds number $90$ around the unstable non-zero steady state;
see~\cite{BehBH17} for details.
The spatial discretization is obtained with Taylor-Hood finite
elements resulting in open-loop systems of the
forms~\cref{eqn:normsys,eqn:nnormsys} described by differential-algebraic
equations, i.e., the $E$ matrices are singular.
The model matrices have been obtained in differently sized
discretizations using the codes from~\cite{BehBH17}.
The resulting dimensions of the systems are given in the cylinder example
column of \Cref{tab:examples}.

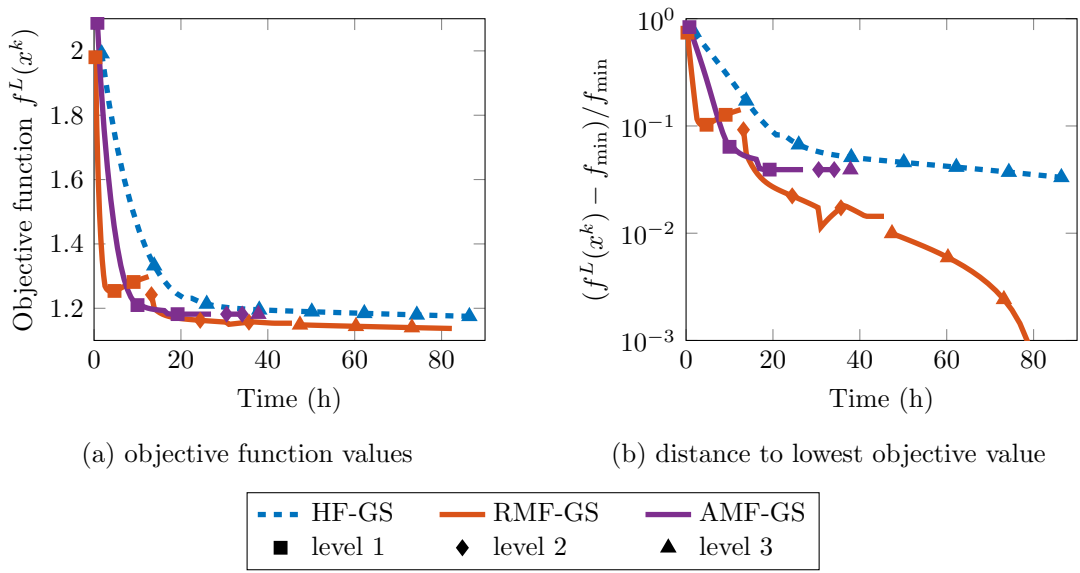
\begin{figure}[t]
  \centering
  \begin{subfigure}[b]{.49\linewidth}
    \centering
  \tikzexternalenable%
  \tikzsetnextfilename{cylinder_norm_objtime}%
  \begin{tikzpicture}[font = \small]
  \begin{axis}[%
    name   = states,
    width  = .7\textwidth,
    height = .18\textheight,
    scale only axis,
    xmin = 0,
    xmax = 90,
    ymin = 1.1,
    ymax = 2.1,
    xminorticks = false,
    yminorticks = false,
    xlabel = {Time (h)},
    ylabel = {Objective function $f^{\maxdisc}(x^{k})$},
    ylabel style   = {yshift = -.3em},
    scaled x ticks = false,
    x tick label style = {/pgf/number format/1000 sep={\,}},
    y tick label style = {/pgf/number format/1000 sep={\,}},
    cycle list name    = cylinderplotlist
  ]
    \pgfplotstableread{graphics/data/cylinder_norm_objtime_gs.dat}\tableDIRECT
    \addplot+[mark repeat = 7] table[x index = 0, y index = 1] {\tableDIRECT};
    
    \pgfplotstableread{graphics/data/cylinder_norm_objtime_rmfgs1.dat}\tableMULTI
    \addplot+[mark repeat = 14] table[x index = 0, y index = 1] {\tableMULTI};
    
    \pgfplotstableread{graphics/data/cylinder_norm_objtime_rmfgs2.dat}\tableMULTI
    \addplot+[mark repeat = 17] table[x index = 0, y index = 1] {\tableMULTI};
    
    \pgfplotstableread{graphics/data/cylinder_norm_objtime_rmfgs3.dat}\tableMULTI
    \addplot+[mark repeat = 7] table[x index = 0, y index = 1] {\tableMULTI};
    
    \pgfplotstableread{graphics/data/cylinder_norm_objtime_amfgs1.dat}\tableMULTI
    \addplot+[mark repeat = 12] table[x index = 0, y index = 1] {\tableMULTI};
    
    \pgfplotstableread{graphics/data/cylinder_norm_objtime_amfgs2.dat}\tableMULTI
    \addplot+[mark repeat = 1] table[x index = 0, y index = 1] {\tableMULTI};
    
    \pgfplotstableread{graphics/data/cylinder_norm_objtime_amfgs3.dat}\tableMULTI
    \addplot+[mark repeat = 1] table[x index = 0, y index = 1] {\tableMULTI};
  \end{axis}
\end{tikzpicture}%
  \tikzexternaldisable%

    \caption{objective function values}
    \label{fig:cylinder_norm_objtime}
  \end{subfigure}%
  \hfill%
  \begin{subfigure}[b]{.49\linewidth}
    \centering
  \tikzexternalenable%
  \tikzsetnextfilename{cylinder_norm_acctime}%
  \begin{tikzpicture}[font = \small]
  \begin{semilogyaxis}[%
    name   = states,
    width  = .7\textwidth,
    height = .18\textheight,
    scale only axis,
    xmin = 0,
    xmax = 90,
    ymin = 1e-03,
    ymax = 1e+00,
    xminorticks = false,
    yminorticks = false,
    xlabel = {Time (h)},
    ylabel = {$(f^{\maxdisc}(x^{k}) - f_{\min}) / f_{\min}$},
    ylabel style   = {yshift = -.3em},
    scaled x ticks = false,
    x tick label style = {/pgf/number format/1000 sep={\,}},
    y tick label style = {/pgf/number format/1000 sep={\,}},
    cycle list name    = cylinderplotlist
  ]
    \pgfplotstableread{graphics/data/cylinder_norm_acctime_gs.dat}\tableDIRECT
    \addplot+[mark repeat = 7] table[x index = 0, y index = 1] {\tableDIRECT};
    
    \pgfplotstableread{graphics/data/cylinder_norm_acctime_rmfgs1.dat}\tableMULTI
    \addplot+[mark repeat = 14] table[x index = 0, y index = 1] {\tableMULTI};
    
    \pgfplotstableread{graphics/data/cylinder_norm_acctime_rmfgs2.dat}\tableMULTI
    \addplot+[mark repeat = 17] table[x index = 0, y index = 1] {\tableMULTI};
    
    \pgfplotstableread{graphics/data/cylinder_norm_acctime_rmfgs3.dat}\tableMULTI
    \addplot+[mark repeat = 7] table[x index = 0, y index = 1] {\tableMULTI};
    
    \pgfplotstableread{graphics/data/cylinder_norm_acctime_amfgs1.dat}\tableMULTI
    \addplot+[mark repeat = 12] table[x index = 0, y index = 1] {\tableMULTI};
    
    \pgfplotstableread{graphics/data/cylinder_norm_acctime_amfgs2.dat}\tableMULTI
    \addplot+[mark repeat = 1] table[x index = 0, y index = 1] {\tableMULTI};
    
    \pgfplotstableread{graphics/data/cylinder_norm_acctime_amfgs3.dat}\tableMULTI
    \addplot+[mark repeat = 1] table[x index = 0, y index = 1] {\tableMULTI};
  \end{semilogyaxis}
\end{tikzpicture}%
  \tikzexternaldisable%

    \caption{distance to lowest objective value}
    \label{fig:cylinder_norm_acctime}
  \end{subfigure}
  
  \vspace{.5\baselineskip}
  \tikzexternalenable%
  \tikzsetnextfilename{cylinder_legend}%
  \begin{tikzpicture}[font = \small]
  \begin{axis}[%
    hide axis,
    width  = .7\textwidth,
    height = .1\textheight,
    scale only axis,
    xmin = 0,
    xmax = 1,
    ymin = 0,
    ymax = 1,
    legend columns = 3,
    legend cell align = {left},
    legend style = {
      at     = {(0,0)},
      anchor = center,
      /tikz/every even column/.append style = {column sep = 0.5cm}}
  ]
    
    \addlegendimage{gs}
    \addlegendentry{\hfgs{}}
    
    \addlegendimage{rmfgs}
    \addlegendentry{\rmfgs{}}
    
    \addlegendimage{amfgs}
    \addlegendentry{\amfgs{}}
    
    \pgfplotsset{cycle list name = cylinderplotlist, cycle list shift = 1}
    \pgfplotsinvokeforeach{1, 2, 3}{
      \addplot+[only marks, black] coordinates {(0,0)};
    }
    \addlegendentry{level $1$}
    \addlegendentry{level $2$}
    \addlegendentry{level $3$}
  \end{axis}
\end{tikzpicture}%
  \tikzexternaldisable%

  \caption{Cylinder example with formulation~\cref{eqn:normsys}: 
    \amfgs{} requires less than half of the runtime time of \hfgs{} to
    converge but it converges to a different objective function value, 
    higher by a factor of about
    $1.0058$.
    \rmfgs{} finds the lowest final objective function value of all three
    methods.}
  \label{fig:cylinder_norm}
\end{figure}

\begin{table}[t]
  \caption{Cylinder example in formulation~\cref{eqn:normsys}:
    The table reports the wall-clock time of the computations, the
    number of iterations taken versus the maximum allowed number and
    the objective function values corresponding to the low-fidelity models
    (in case of \rmfgs{}) and high-fidelity models.}
  \label{tab:cylinder_norm}
  \centering
  {\renewcommand{\arraystretch}{1.25}%
  \settowidth{\maxsteps}{$110$}
  \setlength{\tabcolsep}{.5em}
  \begin{tabular}{lrrrrr}
    \hline
    & &
      Time (h) &
      Iters./Max.\ Iters. &
      $f^{\ell}(x^{k_{\ell}})$ &
      $f^{\maxdisc}(x^{k_{\ell}})$ \\
    \hline\noalign{\medskip}
    \hfgs{} &
      & $86.390$
      & $50$~/~\makebox[\maxsteps][r]{$50$}
      & ---
      & $1.174923$ \\
    \noalign{\medskip}\hline\noalign{\medskip}
    \rmfgs{} & level $1$
      & $12.552$
      & $40$~/~\makebox[\maxsteps][r]{$40$}
      & $1.399568$
      & $1.298534$ \\
    & level $2$
      & $33.007$
      & $50$~/~\makebox[\maxsteps][r]{$50$}
      & $1.152272$
      & $1.153551$ \\
    & level 3
      & $36.802$
      & $20$~/~\makebox[\maxsteps][r]{$20$}
      & ---
      & $1.137206$ \\
    \noalign{\smallskip}\cline{2-6}\noalign{\smallskip}
    & 
      & $82.360$
      & $110$~/~\makebox[\maxsteps][r]{$110$}
      & ---
      & $1.137206$ \\
    \noalign{\smallskip}\hline\noalign{\medskip}
    \amfgs{} & level $1$
      & $26.924$
      & $35$~/~\makebox[\maxsteps][r]{$40$}
      & ---
      & $1.181741$ \\
    & level $2$
      & $7.1769$
      & $2$~/~\makebox[\maxsteps][r]{$50$}
      & ---
      & $1.181741$ \\
    & level $3$
      & $3.7238$
      & $1$~/~\makebox[\maxsteps][r]{$20$}
      & ---
      & $1.181741$ \\
    \noalign{\smallskip}\cline{2-6}\noalign{\smallskip}
    & 
      & $37.852$
      & $38$~/~\makebox[\maxsteps][r]{$110$}
      & ---
      & $1.181741$ \\
    \noalign{\smallskip}\hline\noalign{\smallskip}
  \end{tabular}}
\end{table}

We first consider the formulation~\cref{eqn:normsys}.
The results of the computations can be found in
\Cref{fig:cylinder_norm,tab:cylinder_norm}.
The visible gaps in the lines of \rmfgs{} and \amfgs{} in
\Cref{fig:cylinder_norm} result from the amount of computation time needed to
switch between levels and to perform the first optimization step on the next
level.
The \rmfgs{} method provides the lowest final objective function value of
all methods within about the same runtime as \hfgs{}.
\amfgs{} converges in less than half of the runtime than that of
\rmfgs{} and \hfgs{} but to a different objective function value
than the one found by the other 2 methods, higher by a factor of about $1.0058$.
\amfgs{} finds a good approximation to a stationary point 
already for $\ell = 1$, which cannot be improved further by taking more accurate
gradient sampling steps.
\Cref{tab:cylinder_norm} shows exactly this with its reported numbers of
iterations since for $\ell = 2$, only two steps are performed (one to decrease
the target tolerances of the algorithm and one to verify that no better point
can be found) and only one step for $\ell = 3$, which just confirms that the
approximate stationary point cannot be improved using the given target
tolerances.
However, this point appears to be approximating a local minimizer, as
is indicated by the other two methods obtaining smaller objective function
values.
An interesting point to observe here that we did not see earlier is
that for \rmfgs{}, the high-fidelity objective function value
$f^{\maxdisc}(x^{k})$ is not monotonically decreasing as $k$ increases.
Particularly between $5$ and $15$\,h, the high-fidelity function value
$f^\maxdisc$ increases.
This indicates a mismatch in the approximation of the high-fidelity
model by the low-fidelity model.
Such convergence behavior cannot occur for \amfgs{}, which directly
optimizes the high-fidelity objective function $f^{\maxdisc}$.
Indeed, in the region between $10$ and $15$\,h, the objective function values
obtained by \amfgs{} are smaller than for \rmfgs{} and \hfgs{}.
However, when the discretization is refined,
\rmfgs{} overtakes \amfgs{} and eventually obtains a significantly better
result.
As previously, all three methods needed only a single gradient sampling step to
stabilize the initial controller.

\begin{figure}[t]
  \centering
  \begin{subfigure}[b]{.49\linewidth}
    \centering
  \tikzexternalenable%
  \tikzsetnextfilename{cylinder_objtime}%
  \begin{tikzpicture}[font = \small]
  \begin{axis}[%
    name   = states,
    width  = .7\textwidth,
    height = .18\textheight,
    scale only axis,
    xmin = 0,
    xmax = 45,
    ymin = 0.65,
    ymax = 1.25,
    xminorticks = false,
    yminorticks = false,
    xlabel = {Time (h)},
    ylabel = {Objective function $f^{\maxdisc}(x^{k})$},
    ylabel style   = {yshift = -.3em},
    scaled x ticks = false,
    x tick label style = {/pgf/number format/1000 sep={\,}},
    y tick label style = {/pgf/number format/1000 sep={\,}},
    cycle list name    = cylinderplotlist
  ]
    \pgfplotstableread{graphics/data/cylinder_objtime_gs.dat}\tableDIRECT
    \addplot+[mark repeat = 7] table[x index = 0, y index = 1] {\tableDIRECT};
    
    \pgfplotstableread{graphics/data/cylinder_objtime_rmfgs1.dat}\tableMULTI
    \addplot+[mark repeat = 14] table[x index = 0, y index = 1] {\tableMULTI};
    
    \pgfplotstableread{graphics/data/cylinder_objtime_rmfgs2.dat}\tableMULTI
    \addplot+[mark repeat = 17] table[x index = 0, y index = 1] {\tableMULTI};
    
    \pgfplotstableread{graphics/data/cylinder_objtime_rmfgs3.dat}\tableMULTI
    \addplot+[mark repeat = 7] table[x index = 0, y index = 1] {\tableMULTI};
    
    \pgfplotstableread{graphics/data/cylinder_objtime_amfgs1.dat}\tableMULTI
    \addplot+[mark repeat = 14] table[x index = 0, y index = 1] {\tableMULTI};
    
    \pgfplotstableread{graphics/data/cylinder_objtime_amfgs2.dat}\tableMULTI
    \addplot+[mark repeat = 17] table[x index = 0, y index = 1] {\tableMULTI};
    
    \pgfplotstableread{graphics/data/cylinder_objtime_amfgs3.dat}\tableMULTI
    \addplot+[mark repeat = 7] table[x index = 0, y index = 1] {\tableMULTI};
  \end{axis}
\end{tikzpicture}%
  \tikzexternaldisable%

    \caption{objective function values}
    \label{fig:cylinder_objtime}
  \end{subfigure}%
  \hfill%
  \begin{subfigure}[b]{.49\linewidth}
    \centering
  \tikzexternalenable%
  \tikzsetnextfilename{cylinder_acctime}%
  \begin{tikzpicture}[font = \small]
  \begin{semilogyaxis}[%
    name   = states,
    width  = .7\textwidth,
    height = .18\textheight,
    scale only axis,
    xmin = 0,
    xmax = 45,
    ymin = 1e-04,
    ymax = 1e+00,
    xminorticks = false,
    yminorticks = false,
    xlabel = {Time (h)},
    ylabel = {$(f^{\maxdisc}(x^{k}) - f_{\min}) / f_{\min}$},
    ylabel style   = {yshift = -.3em},
    scaled x ticks = false,
    x tick label style = {/pgf/number format/1000 sep={\,}},
    y tick label style = {/pgf/number format/1000 sep={\,}},
    cycle list name    = cylinderplotlist
  ]
  
    \pgfplotstableread{graphics/data/cylinder_acctime_gs.dat}\tableDIRECT
    \addplot+[mark repeat = 7] table[x index = 0, y index = 1] {\tableDIRECT};
    
    \pgfplotstableread{graphics/data/cylinder_acctime_rmfgs1.dat}\tableMULTI
    \addplot+[mark repeat = 14] table[x index = 0, y index = 1] {\tableMULTI};
    
    \pgfplotstableread{graphics/data/cylinder_acctime_rmfgs2.dat}\tableMULTI
    \addplot+[mark repeat = 17] table[x index = 0, y index = 1] {\tableMULTI};
    
    \pgfplotstableread{graphics/data/cylinder_acctime_rmfgs3.dat}\tableMULTI
    \addplot+[mark repeat = 7] table[x index = 0, y index = 1] {\tableMULTI};
    
    \pgfplotstableread{graphics/data/cylinder_acctime_amfgs1.dat}\tableMULTI
    \addplot+[mark repeat = 14] table[x index = 0, y index = 1] {\tableMULTI};
    
    \pgfplotstableread{graphics/data/cylinder_acctime_amfgs2.dat}\tableMULTI
    \addplot+[mark repeat = 17] table[x index = 0, y index = 1] {\tableMULTI};
    
    \pgfplotstableread{graphics/data/cylinder_acctime_amfgs3.dat}\tableMULTI
    \addplot+[mark repeat = 7] table[x index = 0, y index = 1] {\tableMULTI};
  \end{semilogyaxis}
\end{tikzpicture}%
  \tikzexternaldisable%

    \caption{distance to lowest objective value}
    \label{fig:cylinder_acctime}
  \end{subfigure}
  
  \vspace{.5\baselineskip}
  \tikzexternalenable%
  \tikzsetnextfilename{cylinder_legend}%
  \begin{tikzpicture}[font = \small]
  \begin{axis}[%
    hide axis,
    width  = .7\textwidth,
    height = .1\textheight,
    scale only axis,
    xmin = 0,
    xmax = 1,
    ymin = 0,
    ymax = 1,
    legend columns = 3,
    legend cell align = {left},
    legend style = {
      at     = {(0,0)},
      anchor = center,
      /tikz/every even column/.append style = {column sep = 0.5cm}}
  ]
    
    \addlegendimage{gs}
    \addlegendentry{\hfgs{}}
    
    \addlegendimage{rmfgs}
    \addlegendentry{\rmfgs{}}
    
    \addlegendimage{amfgs}
    \addlegendentry{\amfgs{}}
    
    \pgfplotsset{cycle list name = cylinderplotlist, cycle list shift = 1}
    \pgfplotsinvokeforeach{1, 2, 3}{
      \addplot+[only marks, black] coordinates {(0,0)};
    }
    \addlegendentry{level $1$}
    \addlegendentry{level $2$}
    \addlegendentry{level $3$}
  \end{axis}
\end{tikzpicture}%
  \tikzexternaldisable%

  \caption{Cylinder example with formulation~\cref{eqn:nnormsys}:
    Both \rmfgs{} and \amfgs{} obtain smaller objective function values
    than \hfgs{} does and in a shorter runtime, corresponding to speedups of $4$
    and $2$, respectively.}
  \label{fig:cylinder}
\end{figure}
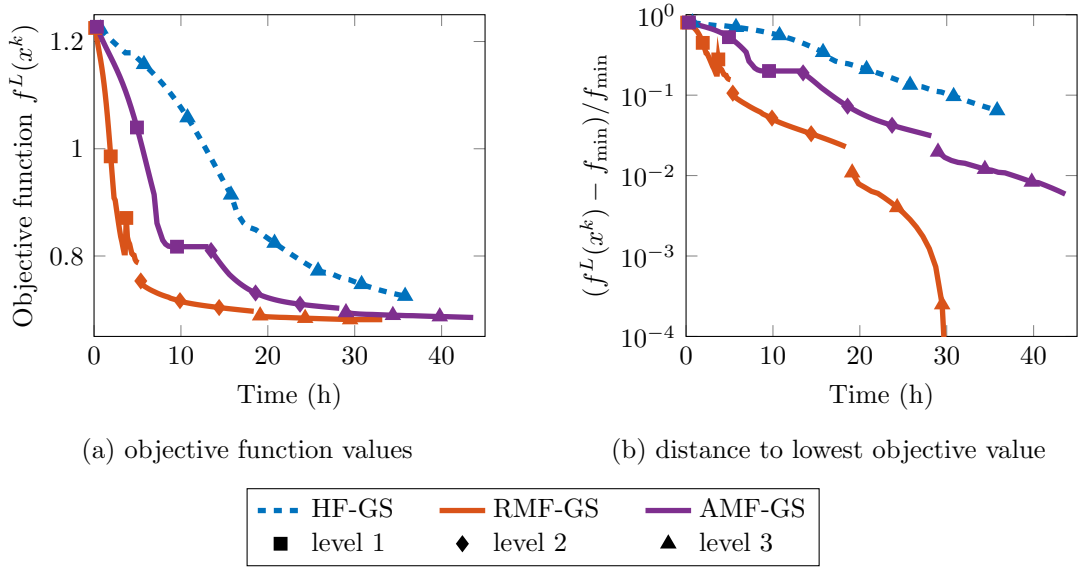

\begin{table}[t]
  \caption{Cylinder example with formulation~\cref{eqn:nnormsys}:
    The table reports the wall-clock time of the computations, the
    number of iterations taken versus the maximum allowed number and
    the objective function values corresponding to the low-fidelity models
    (in case of \rmfgs{}) and high-fidelity models.}
  \label{tab:cylinder}
  \centering
  {\renewcommand{\arraystretch}{1.25}%
  \settowidth{\maxsteps}{$110$}
  \setlength{\tabcolsep}{.5em}
  \begin{tabular}{lrrrrr}
    \hline
    & &
      Time (h) &
      Iters./Max.\ Iters. &
      $f^{\ell}(x^{k_{\ell}})$ &
      $f^{\maxdisc}(x^{k_{\ell}})$ \\
    \hline\noalign{\medskip}
    \hfgs{} &
      & $35.803$
      & $50$~/~\makebox[\maxsteps][r]{$50$}
      & ---
      & $0.725212$ \\
    \noalign{\medskip}\hline\noalign{\medskip}
    \rmfgs{} & level $1$
      & $5.1153$
      & $40$~/~\makebox[\maxsteps][r]{$40$}
      & $0.786143$
      & $0.787230$ \\
    & level $2$
      & $13.262$
      & $50$~/~\makebox[\maxsteps][r]{$50$}
      & $0.696799$
      & $0.696859$ \\
    & level $3$
      & $14.795$
      & $20$~/~\makebox[\maxsteps][r]{$20$}
      & ---
      & $0.681304$ \\
    \noalign{\smallskip}\cline{2-6}\noalign{\smallskip}
    & 
      & $33.172$
      & $110$~/~\makebox[\maxsteps][r]{$110$}
      & ---
      & $0.681304$ \\
    \noalign{\smallskip}\hline\noalign{\medskip}
    \amfgs{} & level $1$
      & $13.169$
      & $40$~/~\makebox[\maxsteps][r]{$40$}
      & ---
      & $0.817351$ \\
    & level $2$
      & $15.055$
      & $50$~/~\makebox[\maxsteps][r]{$50$}
      & ---
      & $0.702572$ \\
    & level $3$
      & $15.415$
      & $20$~/~\makebox[\maxsteps][r]{$20$}
      & ---
      & $0.685316$ \\
    \noalign{\smallskip}\cline{2-6}\noalign{\smallskip}
    & 
      & $43.663$
      & $110$~/~\makebox[\maxsteps][r]{$110$}
      & ---
      & $0.685316$ \\
    \noalign{\smallskip}\hline\noalign{\smallskip}
  \end{tabular}}
\end{table}

Finally, we consider the formulation~\cref{eqn:nnormsys} for the cylinder
example.
The original controls of the benchmark example are modeled to
steer the flow velocities in horizontal and vertical directions behind the
circular obstacle.
We consider only the first half of these controls to introduce disturbances
into the system, which is, for example, the case when control units are
defective.
The second half of the controls remain as given for the design of feedback
controllers.
The results for this example are shown in \Cref{fig:cylinder,tab:cylinder}.
As earlier, \rmfgs{} performs much better than
\amfgs{}, which in turn performs much better than \hfgs{}, obtaining lower
values of $f^{\maxdisc}$ in less runtime.
It requires \amfgs{} $10$\,h more computation time than \rmfgs{} to reach a
value of $f^{\maxdisc}$ that agrees with \rmfgs{} to two digits.
Compared to the final objective function value of \hfgs{}, \amfgs{} performs
around $2$ times faster than \hfgs{} and \rmfgs{} is around $4$ times faster
than \hfgs{}.
For all three methods, only a single gradient sampling step is necessary
to stabilize the initial guess for the controller.

As an alternative to the relatively expensive gradient sampling method, we also
experimented with using the BFGS method, which has proved very effective in
other nonsmooth optimization applications~\cite{CurMO17, LewO13, Ove22}.
However, we found that, particularly for the cylinder example, the behavior of
gradient sampling was more consistent and
reliable, perhaps reflecting its very satisfactory convergence theory, which is
not shared by BFGS.


\section{Conclusions}%
\label{sec:conclusions}

We have introduced two multi-fidelity gradient-sampling approaches for the
robust control of expensive, high-fidelity models that leverage low-cost,
low-fidelity models for speedup.
The numerical experiments demonstrate that speedups of several orders of
magnitude can be achieved compared to a single-fidelity approach that uses
the high-fidelity model alone.
Furthermore, our \rmfgs{} (Restart\-ed Multi-Fidelity Gradient Sampling) method,
which does not access the highest fidelity model until the final phase of the
computation, consistently outperforms our \amfgs{} (Approximate Multi-Fidelity
Gradient Sampling) method, which uses the high-fidelity model throughout the
computation, using lower fidelity gradients in the sampling step.
One might have expected the opposite, since \amfgs{} monotonically reduces
the high-fidelity objective function on the sequence $\{x^{k}\}$.
However, as the cylinder example demonstrated (see \Cref{fig:cylinder_norm}),
even when \rmfgs{} fails to reduce the high-fidelity function on a lower level
of optimization, it can still recover when it continues to the next level of
optimization.
In fact, its robustness seems to reflect its stronger convergence properties.
As explained in \Cref{subsec:theory}, the convergence guarantees of the gradient
sampling algorithm apply at every level of the \rmfgs{} method, while, because
of the approximate gradients used by \amfgs{}, they apply only at the final
level of \amfgs{}, which, in a sense, means that its convergence guarantees are
no stronger than those of \hfgs{}.
One could argue that the consequence of this is that the result of optimization
on one level of \rmfgs{} really does provide a good starting point for
optimization at the next level; the same argument cannot be made for \amfgs{}.

An interesting question that we leave for future work is what convergence
guarantees one might be able to derive for a variant of \rmfgs{} where the
discretization level increases without bound so that it asymptotically
approximates a limit objective function that is computationally intractable.
Such a situation can be found when the dynamical system stems from a
discretization of an underlying partial differential equation and the limit
$\ell \to \infty$ means driving the mesh width to zero to asymptotically
approximate the continuous solution of the partial differential equation and its
corresponding objective function.
Such a setting is considered in the context of uncertainty quantification in,
e.g.,~\cite{CliGSetal11, Gil08, PehGW18}.


\section*{Acknowledgments}%
\addcontentsline{toc}{section}{Acknowledgments}

The authors were partially supported by the National Science Foundation 
under Grant No.~2012250.
The third author was additionally supported by the National
Science Foundation under Grant No.~1901091.
This material is based upon work supported by the National Science
Foundation under Grant No.~DMS-1439786 and by the Simons Foundation Grant
No.~50736 while the first and third author were in residence at the
Institute for Computational and Experimental Research in Mathematics in
Providence, RI, during the ``Model and dimension reduction in uncertain and
dynamic systems'' program.

We would like to thank Tim Mitchell of the Max Planck Institute in Magdeburg,
Germany, who provided a prerelease of ROSTAPACK version~3.0
(now publicly available) for the initial numerical experiments of this work, for
the useful discussions with him about $\Hinf$-norm computations and his valuable
comments on a draft of this manuscript.


\appendix

\section{Gradients of the \texorpdfstring{$\Hinf$}{H-infinity}-norm
  of the closed-loop system}%
\label{app:hinfnorm}
For the use of gradient sampling in $\Hinf$-controller design, the
gradients of the $\Hinf$-norm~\cref{eqn:hinfnorm2} of the closed-loop
system~\cref{eqn:closedsys} with respect to the controller matrices
from~\cref{eqn:K} are needed.
These are well-known in the $\Hinf$-control community, and,
for the case of an identity descriptor matrix in~\cref{eqn:sys}, i.e.,
$E = I_{n}$, they can be found, for example, in~\cite{MilOHetal11}.
We summarize these gradients here for completeness and
include also the case of descriptor matrices as in~\cref{eqn:sys}.
We are concerned with computing the gradients at a given design variable given
by~\cref{eqn:designvar}.
We need to assume that, given these controller variables, the supremum
in~\cref{eqn:linfnorm} is attained only at
one finite point $\omega_{\Hinf}$, with
$\lVert \Gclosed(\i\,\omega_{\Hinf} ) \lVert_{2}
= \lVert \Gclosed \rVert_{\Hinf}$, and that the largest
singular value of $\Gclosed(\i\,\omega_{\Hinf} )$ is simple.
Then the $\Hinf$-norm of the closed-loop system~\cref{eqn:closedsys} is indeed
differentiable and its gradients with respect to the closed-loop
system matrices are given by
\begin{equation} \label{eqn:gradhinfclosed}
  \begin{aligned}
    \nabla_{\Aclosed} \lVert \Gclosed \rVert_{\Hinf} & = Z^{-\herm}
      \Cclosed^{\trans} u v^{\herm} \Bclosed^{\trans} Z^{-\herm}, &
    \nabla_{\Bclosed} \lVert \Gclosed \rVert_{\Hinf} & = Z^{-\herm}
      \Cclosed^{\trans} u v^{\herm}, \\
    \nabla_{\Cclosed} \lVert \Gclosed \rVert_{\Hinf} & = u v^{\herm}
      \Bclosed^{\trans} Z^{-\herm}, &
    \nabla_{\Dclosed} \lVert \Gclosed \rVert_{\Hinf} & = u v^{\herm},
  \end{aligned}
\end{equation}
where $Z = \i\,\omega_{\Hinf}  \Eclosed - \Aclosed$, and $u$ and $v$ are the
right and left singular vectors corresponding to the largest singular value
of $\Gclosed(\i\,\omega_{\Hinf} )$.
Note that the gradient with respect to $\Eclosed$ is not needed since
it does not involve any of the controller matrices, i.e., it contains no
optimization variables for which the gradients need to be evaluated.
However, the matrix $\Eclosed$ plays a role in~\cref{eqn:gradhinfclosed} in
terms of the frequency-dependent matrix pencil $Z$.
Using the chain rule of differentiation we can directly obtain the
requested gradients with respect to the controller matrices
from~\cref{eqn:gradhinfclosed}.
Additionally applying realification to the single terms, since we are
only interested in the design of controllers realized by real-valued
matrices, yields the following results:
\begin{equation} \label{eqn:gradhinf}
  \begin{aligned}
    \nabla_{\AK} \lVert \Gclosed \rVert_{\Hinf} & =
      \real \left( \begin{bmatrix} 0 & I_{\nK} \end{bmatrix}
      \nabla_{\Aclosed} \lVert \Gclosed \rVert_{\Hinf}
      \begin{bmatrix} 0 \\ I_{\nK} \end{bmatrix} \right), \\
    \nabla_{\BK} \lVert \Gclosed \rVert_{\Hinf} & =
      \real \left( \begin{bmatrix} 0 & I_{\nK} \end{bmatrix}
      \nabla_{\Aclosed} \lVert \Gclosed \rVert_{\Hinf}
      \begin{bmatrix} I_{n} \\ 0 \end{bmatrix} C_{2}^{\trans} \right)\\
    & \quad{}+{}
      \real \left( \begin{bmatrix} 0 & I_{\nK} \end{bmatrix}
      \nabla_{\Bclosed} \lVert \Gclosed \rVert_{\Hinf}
      D_{21}^{\trans} \right),\\
    \nabla_{\CK} \lVert \Gclosed \rVert_{\Hinf} & =
      \real \left( B_{2}^{\trans} \begin{bmatrix} I_{n} & 0 \end{bmatrix}
      \nabla_{\Aclosed} \lVert \Gclosed \rVert_{\Hinf}
      \begin{bmatrix} 0 \\ I_{\nK} \end{bmatrix} \right)\\
    & \quad{}+{}
      \real \left( D_{12}^{\trans}
      \nabla_{\Cclosed} \lVert \Gclosed \rVert_{\Hinf}
      \begin{bmatrix} 0 \\ I_{\nK} \end{bmatrix} \right),\\
    \nabla_{\DK} \lVert \Gclosed \rVert_{\Hinf} & =
      \real \left( B_{2}^{\trans} \begin{bmatrix} I_{n} & 0 \end{bmatrix}
      \nabla_{\Aclosed} \lVert \Gclosed \rVert_{\Hinf}
      \begin{bmatrix} I_{n} \\ 0 \end{bmatrix} C_{2}^{\trans} \right)\\
    & \quad{}+{} \real \left(
      B_{2}^{\trans} \begin{bmatrix} I_{n} & 0 \end{bmatrix}
      \nabla_{\Bclosed} \lVert \Gclosed \rVert_{\Hinf}
      D_{21}^{\trans} \right)\\
    & \quad{}+{} \real \left(
      D_{12}^{\trans} \nabla_{\Cclosed} \lVert \Gclosed \rVert_{\Hinf}
      \begin{bmatrix} I_{n} \\ 0 \end{bmatrix} C_{2}^{\trans} \right)\\
    & \quad{}+{} \real \left( D_{12}^{\trans}
      \nabla_{\Dclosed} \lVert \Gclosed \rVert_{\Hinf}
      D_{21}^{\trans} \right).
  \end{aligned}
\end{equation}
Given the $\Hinf$-frequency point $\omega_{\Hinf}$, the gradients
in~\cref{eqn:gradhinf} can be cheaply obtained. This is especially the case
when $\Aclosed$ and $\Eclosed$ are large-scale and sparse by using appropriate
factorizations of the matrix products above.
There have been recent advances in the computation of the $\Linf$-norm
of large-scale sparse systems~\cite{AliBMetal20,BenM18}, which also
yield an efficient approximation of $\omega_{\Hinf}$.


\section{Gradients of the spectral abscissa for initial stabilization}%
\label{app:specabs}
The gradients of~\cref{eqn:specabs} with respect to the controller
matrices~\cref{eqn:K} are well known in the literature for the standard
system case, i.e., $\Eclosed = I_{n + \nK}$; see, for
example,~\cite{BurLO02} and the implementation in~\cite{MilOHetal11}.
Let the design variable be given by~\cref{eqn:designvar}.
We need to assume that the spectral abscissa of
the corresponding matrix pencil $(\Aclosed,\Eclosed)$ is attained at only one
eigenvalue in the closed upper half of the complex plane, say 
$\lambda_{\alpha}$ with $\real(\lambda_{\alpha}) = \alpha(\Aclosed, \Eclosed)$,
and that this eigenvalue is simple.
Then the spectral abscissa is indeed differentiable, with
the gradient, with respect to $\Aclosed$, given by
\begin{equation*}
  \nabla_{\Aclosed} \alpha(\Aclosed, \Eclosed) = w v^{\herm},
\end{equation*}
where $v$ is the right generalized eigenvector of $\lambda_{\alpha}$
and $w$ is the corresponding left eigenvector, normalized
with respect to the inner product with $\Eclosed$, i.e., such that
\begin{equation*}
  w^{\herm} \Eclosed v = 1.
\end{equation*}
Note that we do not need the gradient with respect to $\Eclosed$ since
this matrix does not contain any matrix of the controller~\cref{eqn:K}.
Applying the chain rule and realification of the resulting terms, since
we are only interested in controllers with real-valued matrices, yields the
gradients of interest given by
\begin{align*}
  \nabla_{\AK} \alpha(\Aclosed, \Eclosed) & =
    \real \left( \begin{bmatrix} 0 & I_{\nK} \end{bmatrix}
    \nabla_{\Aclosed} \alpha(\Aclosed, \Eclosed)
    \begin{bmatrix} 0 \\ I_{\nK} \end{bmatrix} \right), \\
  \nabla_{\BK} \alpha(\Aclosed, \Eclosed) & =
    \real \left( \begin{bmatrix} 0 & I_{\nK} \end{bmatrix}
    \nabla_{\Aclosed} \alpha(\Aclosed, \Eclosed)
    \begin{bmatrix} I_{n} \\ 0 \end{bmatrix} C_{2}^{\trans} \right),\\
  \nabla_{\CK} \alpha(\Aclosed, \Eclosed) & =
    \real \left( B_{2}^{\trans} \begin{bmatrix} I_{n} & 0 \end{bmatrix}
    \nabla_{\Aclosed} \alpha(\Aclosed, \Eclosed)
     \begin{bmatrix} 0 \\ I_{\nK} \end{bmatrix} \right),\\
  \nabla_{\DK} \alpha(\Aclosed, \Eclosed) & =
    \real \left( B_{2}^{\trans} \begin{bmatrix} I_{n} & 0 \end{bmatrix}
    \nabla_{\Aclosed} \alpha(\Aclosed, \Eclosed)
    \begin{bmatrix} I_{n} \\ 0 \end{bmatrix} C_{2}^{\trans} \right).
\end{align*}
The right-most eigenvalues and eigenvectors of large-scale sparse matrix
pencils can be efficiently computed using an Arnoldi or Krylov-Schur method
with the shift-and-invert operator and a suitable shift $\sigma$ with a
real part larger than or close to $\alpha(\Aclosed, \Eclosed)$;
see, e.g.,~\cite{GolV13, Ste01}.
The shift $\sigma$ can be efficiently updated during an optimization
approach using the previous computations of $\alpha(\Aclosed, \Eclosed)$.
In our numerical experiments, we use the \texttt{eigs} function from
MATLAB, which in its latest version implements the Krylov-Schur
algorithm~\cite{Ste01}.


\addcontentsline{toc}{section}{References}
\bibliographystyle{plainurl}
\bibliography{bibtex/myref}
  
\end{document}